\documentclass[12pt]{article}
\usepackage{amssymb}
\usepackage{amsfonts}
\usepackage{latexsym}
\usepackage{amsmath}
\usepackage{amsthm}
\def\eb{\begin{eqnarray}}
\def\ee{\end{eqnarray}}
\def\ds{\displaystyle}
\def\non{\nonumber}
\def\wt{\widetilde}

\def\mc{\mathcal}
\def\bp{\begin{proof}[{\bf Proof}]}
\def\g{\gamma}
\def\l{\lambda}
\def\mf{\mathfrak}
\def\D{\mc{D}}

\newtheorem{L1}{Lemma}
\newtheorem{T1}{Theorem}
\newtheorem{P1}{Proposition}
\def\appendix#1{
\addtocounter{section}{1}
\renewcommand{\thesection}{\Alph{section}}
\section*{Appendix \thesection\protect\indent #1}
}

\begin{document}

\begin{titlepage}



\vspace{5 mm}

\begin{center}
{\Large\bf Commuting quantum traces for quadratic algebras}

\vspace{10 mm}

{\bf Zolt\'an Nagy$^{a,}$\footnote{e-mail: nagy@ptm.u-cergy.fr}, Jean Avan$^{a,}$\footnote
{e-mail: avan@ptm.u-cergy.fr}, Anastasia Doikou$^{b,}$\footnote{e-mail: doikou@lapp.in2p3.fr}, 
Genevi\`eve Rollet$^{a,}$\footnote{e-mail: rollet@ptm.u-cergy.fr}}

\vspace{15mm}

{\it $^{a}$Laboratoire de Physique Th\'eorique et Mod\'elisation\\
     Universit\'e de Cergy-Pontoise (CNRS UMR 8089), 5 mail Gay-Lussac, Neuville-sur-Oise,\\
     F-95031 Cergy-Pontoise Cedex, France}
\vspace{5mm}

{\it $^{b}$ Laboratoire d'Annecy-Le-Vieux de Physique Th\'eorique, \\
LAPTH (CNRS UMR 5108), B.P. 110, Annecy-Le-Vieux, F-74941, France}	
\end{center}

\vspace{18mm}

\begin{abstract}
\noindent Consistent tensor products on auxiliary spaces, hereafter denoted ``fusion procedures'', and commuting 
transfer matrices are defined
 for general quadratic algebras, non-dynamical and dynamical, inspired by results
on reflection algebras. Applications of these procedures then yield integer-indexed families of commuting Hamiltonians.

\end{abstract}

\vfill

\end{titlepage}
\setcounter{footnote}{0}
\section{Introduction}

A procedure to construct commuting quantum traces for a particular form of quadratic exchange algebras, known as reflection algebra
\cite{cherednik},
was recently developed in \cite{AD}, building on the pioneering work in \cite{FrMa91a}. We recall that it entails
three different steps: construction of the quadratic exchange algebra itself, and its so-called
``dual'' (this notion will be clarified soon); construction of realizations of the exchange algebra
and its dual on consistent tensor products of the initial auxiliary space (which we will denote here as ``fusion'' procedure)
while keeping a single ``quantum'' Hilbert space on which all operators are assumed to act;
combination of these realizations into traces over the tensorized auxiliary spaces, yielding commuting operators acting on the 
original quantum space, labeled by the integer set of tensorial powers of the auxiliary space.

We immediately insist that this procedure is distinct of, and in a sense complements, the familiar construction of transfer matrices 
by tensoring over distinct {\it quantum} spaces (using an appropriate 
comodule structure of the quantum algebra) while keeping a {\it single} common auxiliary space ; 
the trace is then taken over the auxiliary space to
yield a generating functional of commuting operators\cite{Fad}. In the case when there exists 
a universal formulation of the algebra as a bialgebra with a coproduct structure, both constructions stem from two separate applications of 
this coproduct. 
However, the resulting operators are quite distinct: the trace of the monodromy matrix yields commuting operators
acting on a tensor product of Hilbert spaces (as in e.g. the case of spin chains); the trace of the fused auxiliary matrix yields
operators acting on \emph{one} single Hilbert space. These can be shown in some particular cases to realize the quantum 
analogue of the classical Poisson-commuting traces of powers of the classical Lax-matrix $Tr(L^n)$ (see \cite{AD,ABB,Maillet}).
This is the reason for our phrasing of ``quantum traces'' actually borrowed from \cite{Ma85}. In addition it must be emphasized 
that the procedure itself, combining a construction of a ``dual'' algebra and the establishing of exact fusion formulas, yields
very interesting results on the quadratic exchange algebra itself, and  its possible identification as a coalgebra (e.g. Hopf 
or quasi-Hopf). As we will later comment, it also plays a central role in the (similarly named) Mezincescu-Nepomechie fusion 
constructions for spin chains \cite{MeNe, Annecy}. 

A word of caution is in order. Throughout the paper, we use the term ``fusion'' in a restrictive sense, insofar as we only consider
the possibility of acting on auxiliary spaces. The general fusion procedure itself has been applied also to 
the quantum spaces, yielding e.g. higher spin interactions \cite{KiRe} or multiparticle bound states $S$-matrices.

Our purpose here is to fully describe the quantum trace procedure for three types of general quadratic algebras. The first
one is the quantum non-dynamical quadratic exchange algebra introduced in \cite{FrMa91a}.
The second one was formulated in \cite{NAR} as a dynamical version of the quadratic exchange algebras in \cite{FrMa91a} with
particular zero-weight conditions. It will be denoted ``semi-dynamical'' here, for reasons to be explicited later.
The third one (similarly denoted here as ``fully dynamical'') was first built in \cite{Ra} for the 
$\mathfrak{sl}(2)$ case, and extended to the $\mathfrak{sl}_n$ case in \cite{MPLA}, albeit with particular restrictions on
the coefficient matrices. The zero-weight conditions are different; the algebra structure itself mimicks the reflection algebra introduced
by Cherednik et Sklyanin in \cite{cherednik}; a comodule structure was identified and a universal structure was proposed
in \cite{KM}. We will here briefly comment on the differences 
between the quantum traces built in both dynamical cases.

\section{Non-dynamical quadratic algebras}

These algebras were recognized \cite{cherednik,Ma85} as generalizations of the usual 
$R$-matrix and quantum group structure, leading to non skew symmetrical $r$-matrices in the 
quasiclassical limit. 
 
They are characterised by the following exchange relations.  
\eb
A_{12} \ T_{1} \ B_{12} \ T_{2} \ = \ T_{2} \ C_{12} \ T_{1} \ D_{12} \label{QQA}
\ee
where, as usual, the quantum generators sit in the matrix entries of $T$.
Let us recall some examples of this structure.
\begin{itemize}
\item The Yangian and quantum group structures where
$A=D ,  B=C=
\mathbf{1}$
\item Donin-Kulish-Mudrov (DKM)  reflection algebra without spectral parameters \cite{DKM}.
$A=C  ,  B=D=A^{\pi}$, where $(\phantom{A})^{\pi}$ denotes the permutation of auxiliary spaces: $(A^{\pi})_{12}=A_{21}$.
\item
Kulish-Sklyanin type reflection algebra containing spectral parameters \cite{AD,Sk}:
$A=R_{12}^{-}  ,  B=R_{21}^{+}  ,  C=R_{12}^{+}  ,  D=R_{21}^{-}$ ($\pm$ signs refer to the relative signs of spectral parameters in the
$R$-matrix).
\end{itemize}

In \cite{FrMa91a,FrMa91b} consistency relations involving the structure matrices were derived and it was found that
they had the form of cubic relations on the matrices $A,B,C,D$. 
\eb
A_{12} \ A_{13} \ A_{23} & = & A_{23} \ A_{13} \ A_{12}\label{YBA}\\
A_{12} \ C_{13} \ C_{23} & =& C_{23} \ C_{13} \ A_{12} \label{YBAC}\\
D_{12} \ D_{13} \ D_{23} &=& D_{23} \ D_{13} \ D_{12} \label{YBD}\\
D_{12} \ B_{13} \ B_{23} &=& B_{23} \ B_{13} \ D_{12} \label{YBDB}
\ee
We can see that $A$ and $D$ obey the usual YB-equations whereas $C$ and $B$ are their respective
representations.

Furthermore, generalized unitarity conditions can be derived from self-consistency of (\ref{QQA}) under exchange of 
spaces 1 and 2 which imposes:
\eb\label{unitar}
A_{12}=\alpha A_{21}^{-1} \qquad D_{12}=\beta D_{21}^{-1} \qquad B_{12}=\gamma C_{21} \qquad \left(\alpha,\beta,\gamma \in \mathbb{C}\right)
\ee
The constants of proportionality have to obey an additional 
constraint: $\alpha \gamma=\beta \gamma^{-1}$. In the sequel, we will restrict ourselves to the simplest
choice of $\alpha=\beta=\gamma=1$. 

Let us also note that although $B_{12}=C_{21}$, for \ae sthetical 
and mnemotechnical reasons we continue
to use $C$ whenever it allows for the more familiar and significant $(12,13,23)$ 
display of indices.


In \cite{FrMa91a} the authors had already introduced an algebra which they called  ``dual'' to (\ref{QQA}). 
This ``dual'' structure is characterised by the following exchange relation.\nopagebreak[4] 
\eb
\left(A^{-1}_{12}\right)^{t_1t_2} \ K_1 \ \left(\left(B^{t_1}_{12}\right)^{-1}\right)^{t_2} \ K_2
 = K_2 \left(\left(C_{12}^{t_2}\right)^{-1}\right)^{t_1} \ K_1 \ 
\left(D_{12}^{t_1t_2}\right)^{-1}\label{dQQA}
\ee
Two respective representations of (\ref{QQA}) and (\ref{dQQA}) (assumed to act on 
different quantum spaces) can 
be combined by means of a \emph{trace} \cite{FrMa91a, MeNe,Sk}  on the common auxiliary space to generate commuting 
quantum operators. 
It is with respect to this trace that equation (\ref{dQQA}) 
can be characterized as the  \emph{dual} of equation (\ref{QQA}). We formulate the conjecture that this is
the trace of a $\ast$-algebra structure on some underlying universal algebra. Some freedom remains as to the actual form 
of the trace and in the sequel we will stick to the choice of $H$ as $Tr_{V}(K^{t}T)$. Here the 
superscript $t$ stands for any antimorphism on the auxiliary space $V$, which satisfies also the trace invariance
 property $Tr(KT)=Tr(K^tT^t)$, for all matrices $K$ and $T$. The actual antimorphism may differ from 
the usual transposition (e.g. by additional conjugation, crossing operation) since the proof 
of commutation uses only 
(see theorem ~\ref{theo5},\ref{theo6} and \ref{theo14}) the antimorphism  and  trace invariance properties (see e.g. the super-transposition
in superalgebras, or the crossing operation in $R$-matrices).
Let us also remark here that it is possible to choose a trace formula where the antimorphism 
acts on the quantum space, as it is the case in \cite{AD}, but we prefer not to do so here.
Our  particular choice is motivated by the fact that transposition on the auxiliary space 
is always defined whereas on the quantum space it is not necessarily straightforward and could require a supplementary
hypothesis on this quantum representation which may not be easily implemented.



The quantum trace formulation for such a non-dynamical algebra stems from the results in \cite{AD,FrMa91a}; it is however
interesting to give a rather detailed derivation of it in the general case, since both dynamical algebras will present
similar features, albeit with crucial modifications in the fusion and trace formulas induced by the dynamical
dependence.

We will describe two fusions (consistent tensor product of auxiliary spaces) of equation (\ref{QQA}) respectively
inspired by \cite{AD} (itself relying on \cite{cherednik}) and \cite{DKM}. While the fusion of the structure 
matrices is uniquely defined in each case, 
the solutions of the fused exchange relations are not. In particular, they can be dressed, 
i.e. multiplied by suitable ``coupling'' factors.
This dressing procedure turns out to be crucial: indeed, when the simplest solutions of the fused exchange
relation are combined in a quantum trace, they decouple, giving rise to
products of lower order hamiltonians. To obtain nontrivial commuting quantities these fused $T$-matrices must 
be dressed.

We will finally show that the two fusion procedures identified in \cite{AD,DKM} are related by a coupling 
matrix $L_M$ and that they generate the
same commuting quantities.

\subsection{First fusion procedure \label{1fejez}}
Let us first start by introducing some convenient notations (see \cite{AD}) for fused matrices.

\eb
A_{MN'}=\prod^{\to}_{i \in M} \prod^{\to}_{j \in N'} A_{ij} & =& A_{11'} A_{12'} \ldots A_{1n'}  \non \\
& & A_{21'} A_{22'} \ldots A_{2n'} \ \ldots \non\\
&& A_{m1'} \ldots A_{mn'}
\ee

where $M=\langle 1,2,\ldots,m\rangle$ and $N'=\langle 1',2', \ldots, n'\rangle$ are ordered sets of labels. The same 
sets with 
reversed ordering are denoted
by $\bar{M}$ and $\bar{N'}$. A set $M$ deprived of its lowest (highest) element is denoted by $M_{0}$ ($M^{0}$).

{\bf Remark.} In many explicit examples we would have to deal only with one single exchange formula (\ref{QQA}) with two isomorphic
auxiliary spaces. However  our derivation also applies to a situation where more general coupled sets of exchange relations
would occur as $A_{ij} T_i B_{ij} T_j = T_j C_{ij} T_i D_{ij}$ with $\{i,j\} \subset \{1,\ldots, m_0<\infty \}$ and generically
$V_i \not\approx V_j$. Such situations will occur whenever a universal structure is identifiable and the auxiliary spaces $V_i$ carry
different representations of the algebra, as in e.g. \cite{DKM}. It is therefore crucial that the order in the index set 
be stipulated.

Similar notations are used for the fusion of the other structure matrices.
The next lemma states that that the structure matrices in (\ref{QQA}) can be fused in a way that respects 
the YB-equations (\ref{YBA})-(\ref{YBDB}).
\begin{L1}\label{lemma1}
Let $A, B, C, D$ be solutions of the 
Yang-Baxter equations (\ref{YBA})-(\ref{YBDB}). Then the following fused Yang-Baxter equations hold:
\eb
A_{M\bar{N}'}A_{M\bar{L}''}A_{N'\bar{L}''} &=& A_{N'\bar{L}''}A_{M\bar{L}''}A_{M\bar{N}'}\label{fYBA}\\
A_{M\bar{N}'}C_{ML''}C_{N'L''} &=& C_{N'L''}C_{ML''}A_{M\bar{N}'}\\
D_{M\bar{N}'}D_{M\bar{L}''}D_{N'\bar{L}''} &=& D_{N'\bar{L}''}D_{M\bar{L}''}D_{M\bar{N}'}\\
D_{M\bar{N}'}B_{ML''}B_{N'L''} &=& B_{N'L''}B_{ML''}D_{M\bar{N}'}\label{fYBDB}
\ee
\end{L1}
\bp simple induction on $\#M+\#N'$.\end{proof}

We now describe a fusion procedure for the algebra characterized by (\ref{QQA}), generalizing the one
introduced in \cite{AD}.

\begin{T1}\label{theo1}
If  $T$ is a solution of 
\eb
A_{12} \ T_{1} \ B_{12} \ T_{2} \ = \ T_{2} \ C_{12} \ T_{1} \ D_{12}
\ee
then 
\eb \label{solu1}
 T_{M} = \prod^{\to}_{i \in M} \bigg( T_{i} \bigg( \prod^{\to}_{\substack{i < j\\ j\in M}} 
B_{ij}\bigg)\bigg)
\ee
verifies the following fused equation:
\eb
A_{M\bar{N}'} \ T_{M} \ B_{MN'} \ T_{N'} \ = \ T_{N'} \ C_{MN'} \ T_{M} \ D_{M\bar{N}'}\label{FusEx}
\ee
\end{T1}
\bp
Induction on the cardinality $n$ of the index sets:
$n=\#M+\#N'$ which repeats and generalizes the steps in \cite{AD}. 
\end{proof}

The solution $T_M$ obtained above can be dressed, i.e. can be multiplied from the left and 
the right by suitable factors. 
\begin{P1}\label{prop1}
Let $T_M$ be a solution of the fused exchange relation. Then $Q_M T _M S_M$ is also a solution of the fused exchange relation
provided $Q_M$ and $S_M$ verify:

\eb
\left[Q_M,A_{M\bar{N}'}\right]=\left[Q_{N'},A_{M\bar{N}'}\right]=\left[Q_{N'},B_{MN'}\right]=\left[Q_M,C_{MN'}\right]=0 \label{Dress} \\
\left[S_M,D_{M\bar{N}'}\right]=\left[S_{N'},D_{M\bar{N}'}\right]=\left[S_{N'},C_{MN'}\right]=\left[S_M,B_{MN'}\right]=0 \non
\ee

A particular solution of these constraints is provided by: 
\eb\label{dress1}
Q_{M}=\check{A}_{12}\check{A}_{23}\ldots \check{A}_{m-1,m}\\
S_{M}=\check{D}_{12}\check{D}_{23}\ldots \check{D}_{m-1,m}\non
\ee
where $\check{A}_{12}=P_{12}A_{12}, \ \ldots$, $P_{12}$ being the permutation exchanging two auxiliary spaces.
\end{P1}
\bp again by induction on the cardinality of the index sets. In the induction step we use the decomposition: $Q_{N'}B_{MN'}=
\check{A}_{12}\ldots \check{A}_{n'-1,n'}B_{M,{N'}^{00}}B_{M,n'-1}B_{M,n'}$, for example.\end{proof}

The fusion procedure can be repeated for the dual exchange relation as follows.
\begin{T1}
If $K$ is a solution of the dual exchange relation:
\eb \label{dFusEx}
\left(A^{-1}_{12}\right)^{t_1t_2} \ K_1 \ \left(\left(B^{t_1}_{12}\right)^{-1}\right)^{t_2} \ K_2
 = K_2 \left(\left(C_{12}^{t_2}\right)^{-1}\right)^{t_1} \ K_1 \ \left(D_{12}^{t_1t_2}\right)^{-1}
\ee
then 
\eb\label{dsolu1}
K_{M}= \prod^{\to}_{i \in M} \bigg( K_{i} \bigg( \prod^{\to}_{\substack{i < j\\j \in M}} \left(\left(B^{t_1}_{ij}\right)^{-1}\right)^{t_2}
\bigg)\bigg)
\ee
is a solution of the dual fused equation
\eb\label{dNQAf} 
(A^{-1}_{M\bar{N}'})^{t_Mt_{N'}} K_M  \left((B^{t_M}_{MN'})^{-1}\right)^{t_{N'}}  K_{N'}
 = K_{N'} \left((C_{MN'}^{t_{N'}})^{-1}\right)^{t_M}  K_M  (D_{M\bar{N}'}^{t_Mt_{N'}})^{-1}
\ee
\end{T1}
\bp similar to that of Theorem \ref{theo1}. 
Note that the dual structure matrices obey a set of appropriate YB-equations, isomorphic to (\ref{fYBA})-(\ref{fYBDB}), for instance
\eb(A^{-1}_{M\bar{N}'})^{t_Mt_{N'}}
(A^{-1}_{M\bar{L}''})^{t_Mt_{L''}} (A^{-1}_{N\bar{L}''})^{t_Nt_{L''}} = 
(A^{-1}_{N\bar{L}''})^{t_Nt_{L''}}
(A^{-1}_{M\bar{L}''})^{t_Mt_{L''}} (A^{-1}_{M\bar{N}'})^{t_Mt_{N'}} \non \\ 
 \nopagebreak[4]
\ee
\end{proof}
A similar dual dressing procedure exists:
Any dressing of a solution of (\ref{dFusEx}) should obey the commutativity constraints 
\eb
\left[Q'_M,(A^{-1}_{M\bar{N}'})^{t_Mt_{N'}}\right]=\left[Q'_{N'},(A^{-1}_{M\bar{N}'})^{t_Mt_{N'}}\right]=
\left[Q'_{N'},((B^{t_M}_{MN'})^{-1}\right]=\\ \non \left[Q'_M,((C_{MN'}^{t_{N'}})^{-1}\right]=  
\left[S'_M,(D_{M\bar{N}'}^{t_Mt_{N'}})^{-1}\right]=\left[S'_{N'},(D_{M\bar{N}'}^{t_Mt_{N'}})^{-1}\right]=\\
\left[S'_{N'},((C_{MN'}^{t_{N'}})^{-1}\right]=\left[S'_M,((B^{t_M}_{MN'})^{-1}\right]=0 \non
\ee
involving fused dual structure matrices. It is easy to check that 
if $Q_M$ and $S_M$ dress solutions of (\ref{FusEx}) then $Q'_M=Q_M^{t_M}$ and $S'_M=S_M^{t_M}$ 
dress solutions of (\ref{dFusEx}).

\subsection{Second fusion procedure}

Results in \cite{DKM} hint that relation (\ref{QQA}) admits another fusion procedure. We will explicitely
link the fusion described in the preceding section to the one inspired by ref. \cite{DKM}. 

The DKM type fusion is characterized by the following fused exchange relation for fused matrices $\wt{T}$ to be
described in the following:
\eb
A_{\bar{M}N'} \ \wt{T}_{M} \ B_{M\bar{N}'} \ \wt{T}_{N'} \ = \ \wt{T}_{N'} \ C_{\bar{M}N'} \ \wt{T}_{M} 
\ D_{M\bar{N}'}
\ee
This equation can actually be obtained from a multiplication of the KS exchange relation (\ref{FusEx}) by suitable 
factors reversing the ordering of indices where it is needed. 
The next lemma specifies this statement.
\begin{L1}\label{lemma2}
Let $T_M$ be a solution of the fused exchange relation (\ref{FusEx}). If
 $L_M$ verifies the following commutation rules  
\eb \label{kernel}
&L_MA_{M\bar{N}'}=A_{\bar{M}\bar{N}'} L_M \, \qquad  \qquad L_{N'}A_{M\bar{N}'} = A_{MN'} L_{N'} &\\
&L_{N'} B_{MN'} =B_{M\bar{N}'} L_{N'} \,  \qquad  \qquad L_{M} C_{MN'}=C_{\bar{M}N'} L_{M}&\non
\ee
then $\wt{T}_M=L_MT_M$ is a solution of the exchange relation
\eb
A_{\bar{M}N'} \ \wt{T}_{M} \ B_{M\bar{N}'} \ \wt{T}_{N'} \ = \ \wt{T}_{N'} \ C_{\bar{M}N'} \ \wt{T}_{M} 
\ D_{M\bar{N}'}\label{FusEx2}
\ee
An example of such an $L_M$ is given by:
\eb\label{kernelEx}
L_M=A_{12}\ldots A_{1m} A_{23}\ldots A_{2m} \ldots A_{m-1,m} = 
\prod^{\rightrightarrows}_{1\le i < j \le m} A_{ij}
\ee
\end{L1}
\bp The first part is straightforward. Example (\ref{kernelEx}) is verified by induction using
$L_{M}=A_{1M_0}L_{M_0}$. For instance, the first relation of (\ref{kernel}) is proved as:
\eb
&&L_MA_{M\bar{N}'}=A_{1M_0} L_{M_0} A_{1\bar{N}'} A_{M_0\bar{N}'} = A_{1M_0}A_{1\bar{N}'} L_{M_0} A_{M_0\bar{N}'} = \non \\
&&A_{1M_0} A_{1\bar{N}'} A_{\bar{M}_0\bar{N}'} L_{M_0} = A_{\bar{M}_0\bar{N}'} A_{1\bar{N}'} A_{1M_0} = 
A_{\bar{M}\bar{N}'} L_{M} \non
\ee
where fused YB-equations are used.\end{proof}

Combined with Theorem \ref{theo1}, this lemma leads to 
\begin{T1}\label{theo3}
If $\wt{T}$ is a solution of 
\eb
A_{12} \ \wt{T}_{1} \ B_{12} \ \wt{T}_{2} \ = \ \wt{T}_{2} \ C_{12} \ \wt{T}_{1} \ D_{12} 
\ee
then
\eb \label{solu2}
\wt{T}_M=\prod^{\to}_{i \in M}\bigg(\prod^{\to}_{\substack{j>i \\ j \in M}}A_{ij} \wt{T}_i \prod^{\gets}_{\substack{j>i \\ j \in M}} B_{ij}\bigg)
\ee
is a solution of 
\eb
A_{\bar{M}N'} \ \wt{T}_{M} \ B_{M\bar{N}'} \ \wt{T}_{N'} \ = \ \wt{T}_{N'} \ C_{\bar{M}N'} \ \wt{T}_{M} 
\ D_{M\bar{N}'}
\ee
\end{T1}
\bp The only property left to check is that the solution $\wt{T}_M$ in (\ref{solu2}) is obtained from
$T_M$ in (\ref{solu1}) by a multiplication by $L_M$ in (\ref{kernelEx}). It is enough to show that 
$\wt{T}_M=A_{1M_0} \wt{T}_1 B_{1\bar{M}_0} \wt{T}_{M_0}$. We only develop the induction step.
\eb
&&L_M T_M= A_{1M_0} L_{M_0}T_1 B_{1M_0} T_{M_0}= A_{1M_0} T_1 L_{M_0} B_{1M_0}T_{M_0}= \\
&&A_{1M_0} T_1 B_{1\bar{M}_0} L_{M_0} T_{M_0} \non
=A_{1M_0} \wt{T}_1 B_{1\bar{M}_0} \wt{T}_{M_0} 
\ee
\end{proof}

The next proposition describes the dressing of the solutions.
\begin{P1}
Let $\wt{T}_M$ be a solution of the DKM-type fused exchange relations. Then $\wt{Q}_M \wt{T}_M \wt{S}_M$ is also 
a solution provided $\wt{Q}_M$ and $\wt{S}_M$ verify
\eb
\left[\wt{Q}_M,A_{\bar{M}N'}\right]=\left[\wt{Q}_{N'},A_{\bar{M}N'}\right]=\left[\wt{Q}_{N'},B_{M\bar{N}'}\right]
=\left[\wt{Q}_M,C_{\bar{M}N'}\right]=0 \label{Dress2} \\
\left[\wt{S}_M,D_{M\bar{N}'}\right]=\left[\wt{S}_{N'},D_{M\bar{N}'}\right]=\left[\wt{S}_{N'},C_{\bar{M}N'}\right]
=\left[\wt{S}_M,B_{M\bar{N}'}\right]=0 \non
\ee
These equations are solved by
\eb
&\wt{Q}_M=L_M  \ Q_M  \ L_M^{-1} \qquad
\wt{S}_M=S_M \non
\ee
where $Q_M$ and $S_M$ dress the solutions of the fused exchange relation (\ref{FusEx}) and $L_M$ 
is a solution of (\ref{kernel}).
\end{P1}
\bp Straightforward.\end{proof}

We saw that $T_M$ and $\wt{T}_M$ were linked by a factor $L_M$. The question arises whether there 
is a similar relation between the corresponding dual exchange algebras and their solutions. 
The relation is established in
\begin{T1}\label{theo4}
Let $K_M$ be a solution of the first fused exchange relation (\ref{FusEx}) and $L_M$ be a solution of (\ref{kernel}).
Then $\wt{K}_M=(L_M^{t_M})^{-1} K_M$ is a solution of the KDM-type dual fused exchange relation:
\eb \label{dFusEx2}
(A^{-1}_{\bar{M}N'})^{t_Mt_{N'}} \wt{K}_M  ((B^{t_M}_{M\bar{N}'})^{-1})^{t_{N'}}  \wt{K}_{N'}
 = \wt{K}_{N'} ((C_{\bar{M}N'}^{t_{N'}})^{-1})^{t_M}  \wt{K}_M  (D_{M\bar{N}'}^{t_Mt_{N'}})^{-1}
\ee
\end{T1}

\bp We first see that (\ref{dFusEx2}) is indeed the dual exchange relation associated with (\ref{FusEx2}). The next step is
to check that $(L_M^{t_M})^{-1}$ obeys the appropriate commutation relations that enable it to transform the fused dual AD type 
algebra (\ref{dFusEx}) into the fused dual DKM-type one (\ref{dFusEx2}). It is
obvious since these equations are the inverse-transposed of (\ref{kernel}).\end{proof}

Dressings of these dual fused solutions are obtained from dressings of (\ref{FusEx2}) by the same operation as for the AD type
fusion, i.e. by transposing.

\subsection{Commuting traces}

In the preceding sections we have derived two distinct fusion procedures both of which allow for building commuting
quantities. In this section we will describe this construction, and show the two different 
quantum traces are identified once the dressing is used.

We first establish:
\begin{T1} \label{theo5}

Let $\mc{T}_{M}$ be a solution of the fused AD-type exchange relation (\ref{FusEx}). $\mc{T}_{M}$ acts on the tensor 
product of the auxiliary spaces labeled by $M$ and on the quantum space $V_q$.

Let $\mc{K}_M$ be a solution of the dual fused AD-type exchange relation (\ref{dFusEx}). $\mc{K}_M$ acts on the tensor 
product of the auxiliary  spaces labeled by $M$ and on the quantum space $V_{q'}$.

The following operators 
\eb
H_M=Tr_M \left(\mc K_M^{t_M} \mc T_M\right) \label{ComHam}
\ee
constitute a family of mutually commuting quantum operators acting on $V_q \otimes V_{q'}$:
\eb
\left[H_M,H_{N'}\right]=0
\ee
\end{T1}

\bp It repeats the steps of \cite{AD,Sk} 
\end{proof}

The proof is independent of the particular fusion procedure so it remains valid for the DKM case too. 
Thus we have

\begin{T1}
\label{theo6}
Let $\wt{\mc{T}}_{M}$ be a solution of the fused DKM-type exchange relation (\ref{FusEx2}). $\wt{\mc{T}}_{M}$ acts on the tensor 
product of the auxiliary spaces labeled by $M$ and on the quantum space $V_q$.

Let $\wt{\mc{K}}_M$ be a solution of the dual fused DKM-type exchange relation (\ref{dFusEx2}). $\wt{\mc{K}}_M$ acts on the tensor 
product of the auxiliary  spaces labeled by $M$ and on the quantum space $V_{q'}$.

The following operators 
\eb
\wt H_M=Tr_M \left(\wt{\mc K}_M^{t_M} \wt {\mc T}_M\right)\label{CommHam}
\ee
constitute a family of mutually commuting quantum operators acting on $V_q \otimes V_{q'}$:
\eb
\left[\wt H_M,\wt H_{N'}\right]=0
\ee
\end{T1}\vspace{5mm}

So far we have two seemingly different sets of commuting quantities obtained from the same 
defining relations (\ref{QQA}) via two distinct fusion procedures. However we will show that 
the operation consisting in dressing and taking the trace smears out this
difference and one is left with only one set of commuting hamiltonians. This is summarized in:

\begin{P1}\label{thesame}
The quantum commuting Hamiltonians obtained from any set of solutions $T_M, \ K_M$ of (\ref{FusEx}), 
(\ref{dNQAf}) are identified with the quantum commuting Hamiltonians obtained from a suitable
set of solutions $\tilde{T}_M, \ \tilde{K}_M$ of (\ref{FusEx2}), (\ref{dFusEx2}). This identification
is implemented by a coupling matrix $L_M$.
\end{P1}
\bp Let $\mc T_M$ be the solution (\ref{solu1}) and $\mc K_M$ the corresponding dual solution (\ref{dsolu1}). 
The results of the multiplication by $L_M$ and $(L_M^{t_M})^{-1}$  are denoted by $\wt{\mc T_M}$ and $\wt{\mc K}_M$. 
We calculate the tilded hamiltonians after dressing and we find that they are equal to the dressed untilded ones.
\eb
&&Tr_M(\wt{\mc K}_M^{t_M} \wt Q_M \wt {\mc T}_M \wt S_M) = Tr_M(\mc K_M^{t_M} L_M^{-1} L_M Q_M L_M^{-1} L_M \mc T_M S_M) = \non \\
&&Tr_M(\mc K_M^{t_M} Q_M \mc T_M S_M) \ 
\ee
\end{proof}

The following propositions justifies the technical relevance of dressings.
\begin{P1}\label{decoup}
Operators built from the solution (\ref{solu1}) decouple as $H_N=Tr_N(K_N^{t_N}T_N)=Tr(K^tT)^{\#N}$.
\end{P1}
\bp By induction using the property $T_N=T_1 B_{1N_0} T_{N_0}$. Let us detail the induction step.
\eb
&&H_N=Tr_N(K_N^{t_N}T_N)=Tr((K_1(B^{t_1}_{1N_0})^{-1t_{N_0}}K_{N_0})^{t_{N}}T_1 B_{1N_0}T_{N_0})= \non \\
&&Tr((K_1(B^{t_1}_{1N_0})^{-1t_{N_0}}K_{N_0})^{t_{N_0}}(T_1B_{1N_0}T_{N_0})^{t_1})= \non \\
&&Tr(K_1K_{N_0}^{t_{N_0}}(B^{t_1}_{1N_0})^{-1}B_{1N_0}^{t_1}T_1^{t_1}T_{N_0})= \non \\
&&Tr(K_1T_1^{t_1})Tr(K_{N_0}^{t_{N_0}}T_{N_0}) \non
\ee \end{proof}

Note that the result in Proposition \ref{thesame} implies that the same goes for the operators built using the second fusion.
Three important remarks are in order here.

\subsubsection*{The use of dressed quantum traces}

Dressed quantum traces yield a priori independent operators. Indeed, the classical limit of a quantum trace 
computed with the particular dressing (\ref{dress1}) in Proposition \ref{prop1} will yield $Tr T^n$ 
instead of $\left(Tr T\right)^n$ (since $A, B, C, D 
\rightarrow \mathbf{1} 
\otimes \mathbf{1}$ but $P_{12} \to P_{12}$ \,!). Quantum traces are directly, in this particular case, (as was already known in the 
context of quantum group structures \cite{Maillet}) quantum analogues of the classical Poisson-commuting power traces $Tr T^n$.
\subsubsection*{The use of undressed quantum traces}
It must on the other hand be emphasized that the decoupling of the undressed fused quantities plays an essential role  in the formulation
of the analytical Bethe ansatz solution of $\mathfrak{sl}(n)$ spin chains (as is seen in \cite{Annecy}) and 
more generally in the formulation of a generalized Mezincescu-Nepomechie procedure for fusion of transfer matrices \cite{MeNe}, in
that it gives a natural construction of products of monodromy matrices such as are required by this formulation.

\subsubsection*{Explicit computation of the dressings}

From a more theoretical point of view, it must be noticed that eqn. (\ref{kernel}), as already discussed for the particular
example treated in \cite{AD}, would appear as a condition obeyed by coproducts of the central elements of a 
(hypothetical) universal algebra, thereby promoting
the dressing matrices $Q$ and $S$ from ``technical auxiliaries'' to get non-trivial traces, to representations of
Casimir elements of the algebra itself\footnote{this was pointed out to us by Daniel Arnaudon}.

A second more technical remark is required here regarding the actual  computation of the quantum traces with the particular 
explicit dressing 
determined in Proposition \ref{prop1}. Difficulties in applying (\ref{CommHam}) with the explicit dressings (\ref{dress1}) may occur
when the auxiliary space $V$ is a loop space $V^{(n)}\otimes \mathbb{C}(z)$ ($n$=finite dimension of the vector space). Indeed,
the permutation of spectral parameters required in formula (\ref{dress1}) is only achieved at 
a formal level by the singular distribution $\delta (z_i/z_j)$ (see \cite{AD} for discussions). Hence the actual 
explicit computations of such quantum traces may entail delicate regularization procedures. However, if one only focuses
on the practical purpose of the quantum trace procedure, which is to build a set of commuting operators, use of higher-power
fused objects as in (\ref{solu1}) and (\ref{dsolu1}) is mostly required when no spectral parameter is present in the 
represented exchange algebra (\ref{QQA}). Otherwise one needs to consider only the first order trace 
$Tr_1 \wt{K}_1(z_1)T_1(z_1)$ and expand it in formal series in $z_1$. If no spectral parameter is available, one can then use
(\ref{solu1}), (\ref{dsolu1}), (\ref{dress1}) and (\ref{CommHam}) to build explicitely without difficulties a priori
independent commuting quantum operators. (For an application to a different algebraic structure see \cite{ABB}).

\subsection{Further example: ``Soliton non-preserving'' boundary conditions: Twisted Yangians.}
We have mentioned in the Introduction several examples
of non-dynamical quadratic exchange algebras.
Another interesting example to which we plan to apply this scheme
is related to the so-called
``soliton non-preserving'' boundary conditions in integrable
lattice models (see \cite{doikou}).
To characterize it 
we will focus on the $\mathfrak{su}(n)$ invariant $R$ matrix given by
\eb R_{12}(\lambda) = \lambda I+i{\cal P}_{12} \ee where ${\cal P}$ is the permutation operator on the tensor product $V_1 \otimes V_2$. 
The $R$ matrix is a solution  of the Yang--Baxter equation \cite{mac,yang,baxter,korepin} and also  satisfies:
 
{\it (i) Unitarity}
\eb
R_{12}(\lambda)\ R_{21}(-\lambda) = \zeta(\lambda) \label{uni1} \ee
where $R_{21}(\lambda) ={\cal P}_{12} R_{12}(\lambda) {\cal P}_{12} = R_{12}^{t_{12}}(\lambda)$ and ${\cal P}$ is the permutation operator. 

{\it (ii) Crossing--unitarity} 
\eb
R_{12}^{t_{1}}(\lambda)\ M_{1}\ R_{12}^{t_{2}}(-\lambda-2i\rho)\ M_{1}^{-1} =\zeta'(\lambda+i\rho) \label{croun}\ee 
$M= V^{t}\ V$, ($M=1$ for the $\mathfrak{su}(n)$ case)  $\ds \rho = {\frac{n}{2}}$ and also
\eb [M_{1}M_{2},\ R_{12}(\lambda)] =0, \ee
\eb \zeta(\lambda) = (\lambda+i)(-\lambda +i), ~~\zeta'(\lambda) = (-\lambda +i\rho)(\lambda+i\rho). \ee
It is interpreted as the scattering matrix \cite{zamo, korepin, faddeev}
 describing the interaction between two solitons --objects that correspond
 to the fundamental representation of $\mathfrak{su}(n)$.

 One may also derive
 the scattering matrix that describes the interaction between a soliton and
 an anti-soliton, which corresponds to the conjugate representation of
 $\mathfrak{su}(n)$. It reads:
\eb R_{\bar 12}(\lambda) =R_{1 \bar 2}(\lambda) = \bar R_{12}(\lambda)= U_{1}\ R_{12}^{t_{2}}(-\lambda -i\rho)\ U_{1}, 
\label{cross} \ee 
and it  can also be written as
\eb \bar R_{12}(\lambda) = (-\lambda -i\rho)I+iQ \ee where $Q$ is a projector onto a one dimensional space, and
where $U$ is a matrix of square $1$. 
Note that for the $su(2)$ case \eb \bar R_{12}(\lambda) =R_{12}(\lambda), \ee which is expected because $su(2)$ is self conjugate.
The $\bar R$  matrix also satisfies the Yang--Baxter equation  and 

{\it (i) Unitarity}
\eb
\bar R_{12}(\lambda)\ \bar R_{2 1}(-\lambda) = \zeta'(\lambda) \label{uni2} \ee 

{\it (ii) Crossing--unitarity} 
\eb
\bar  R_{12}^{t_{1}}(\lambda)\ M_{1}\ \bar R_{12}^{t_{2}}(-\lambda-2i\rho)\ M_{1}^{-1} =\zeta(\lambda).  \label{croun2}\ee 

\subsubsection*{The reflection equation}
The usual reflection equation \cite{cherednik} describes physically
 the reflection of a soliton (fundamental representation of $\mathfrak{su}(n)$) as
 a soliton. The associated quadratic algebra was considered e.g. in \cite{AD}
\begin{equation}
R_{12}(\lambda_{1}-\lambda_{2})\ T_{1}(\lambda_{1})
 R_{21}(\lambda_{1}+\lambda_{2})\ 
T_{2}(\lambda_{2})=
T_{2}(\lambda_{2})\ R_{12}(\lambda_{1}+\lambda_{2})\ 
T_{1}(\lambda_{1})\ R_{21}(\lambda_{1}-\lambda_{2}). \label{re1}
\end{equation}
Considering now the reflection of a soliton as anti-soliton one is similarly
lead to the formulation of another quadratic algebra:

\begin{equation}
R_{12}(\lambda_{1}-\lambda_{2})\ T_{1}(\lambda_{1})\
\bar R_{21}(\lambda_{1}+\lambda_{2})\
T_{2}(\lambda_{2})=
T_{2}(\lambda_{2})\ \bar R_{12}(\lambda_{1}+\lambda_{2})\
T_{1}(\lambda_{1})\ R_{21}(\lambda_{1}-\lambda_{2}). \label{re2}
\end{equation} More
specifically equation (\ref{re2}) is the definition of the so--called
 {\it twisted Yangian}.
Its dual reflection equation is
 obtained essentially by taking its formal transposition:
\eb &&  R_{12}(-\lambda_{1}+\lambda_{2})\ K_{1}^{t_{1}}(\lambda_{1})\ M_{1}^{-1}\
 \bar R_{21}(-\lambda_{1}-\lambda_{2}-2i\rho)\ M_{1}\
K_{2}^{t_{2}}(\lambda_{2}) \non\\=&&
K_{2}^{t_{2}}(\lambda_{2})\ M_{1}\   \bar R_{12}(-\lambda_{1}-\lambda_{2}-2i\rho)\
 M_{1}^{-1}\ K_{ 1}^{t_{1}}(\lambda_{1})\  R_{21}(-\lambda_{1}+\lambda_{2}). \label{re+}
\ee  
This indeed realizes the general quadratic exchange relation (\ref{QQA}),
 (\ref{dQQA}) with the following identifications (using unitarity and crossing symmetries of the $R$-matrix)
\eb A_{12} = R_{12}(\lambda_{1}-\lambda_{2}), ~B_{12}
 =\bar R_{21}(\lambda_{1}+\lambda_{2}),~C_{12}= \bar R_{12}(\lambda_{1}+\lambda_{2}), ~D_{12} =R_{21}(\lambda_{1}-\lambda_{2}) \non 
\ee
\eb && (A_{12}^{-1})^{t_1t_2} = R_{12}(-\lambda_{1}+\lambda_{2}),
 ~ ((B_{12}^{t_1})^{-1})^{t_2} =M_{1}\bar R_{12}(-\lambda_{1}-\lambda_{2}-2i\rho)M_{1}^{-1},
\non\\ && ((C_{12}^{t_2})^{-1})^{t_1}=M_{1}^{-1} \bar R_{21}(-\lambda_{1}-\lambda_{2}-2i\rho)M_{1}, 
~ (D_{12}^{-1})^{t_1t_2} =R_{21}(-\lambda_{1}+\lambda_{2}) \non \ee

Explicit application of the quantum trace procedure to this
particular algebra will be left for further studies.

\section{Quantum traces for semi-dynamical quadratic algebras}

The second type of quadratic exchange relations considered here consists of the dynamical quadratic algebras 
generically described and studied in \cite{NAR} which were first 
exemplified in the context of scalar Ruijsenaars-Schneider models in \cite{ACF}.
Fusion procedures and commuting traces can be built up for these dynamical quadratic algebras following the same overall
procedure as in the non-dynamical case, albeit with crucial, non-tr0ivial differences.

\subsection{The semi-dynamical quadratic algebra}

Let us recall here the basic definitions.
Our starting point is the dynamical quadratic exchange relation:
\eb\label{DQQA}
A_{12}(\lambda)T_1(\lambda)B_{12}(\lambda)T_2(\lambda+\gamma h_1)&=&T_2(\lambda)C_{12}(\lambda)T_1(\lambda+\gamma h_2)D_{12}
\ee
This describes an algebra generated by the matrix entries of $T$.
$A,B,C,D$ are matrices in $End(V\otimes V)$ depending on $\lambda \in \mathfrak{h}^{\ast}$ where $\mathfrak{h}$ is a  commutative
Lie algebra, of dimension $n$, making $V$ a diagonalizable $\mathfrak{h}$-module. 
Introducing coordinates $\lambda_i$ on 
$\mathfrak{h}^{\ast}$ and the dual base $h_i$ on $\mf{h}$ the shift $\lambda +\gamma h$ can  be defined in the following way. 
For any differentiable function $f(\lambda)=f(\{\lambda_i\})$:
\eb
f(\lambda+\gamma h)=e^{\gamma \mc{D}}f(\lambda)e^{-\gamma \mc{D}},
\ee
 where
\eb
\mc{D}=\sum_i h_i \partial_{\lambda_i}
\ee
In the forthcoming calculations $\gamma$ is set to 1 for simplification.
Zero weight conditions are imposed on the first space of $B$ and the second one of $C$; $D$ is of total weight zero. 
\eb
\left[B_{12},h\otimes \mathbf{1}\right]=\left[C_{12},\mathbf{1}\otimes h\right]=\left[D_{12},h \otimes \mathbf{1}+\mathbf{1}\otimes h
\right]=0 \quad (h\in \mathfrak{h})
\ee
These particular conditions, together with the absence of dynamical shift in two out of four $T$ matrices in (\ref{DQQA}), lead
us to denote this structure as ``semi-dynamical''. 
We will restrict ourselves from now on to the case where $V$ is of dimension $n$: the basis of $V$ and the generators of $\mf{h}$,
can then be chosen so that one identifies: 
 $h_i=E_{ii}$ (diagonal basis elements of $\mf{gl}(n)$,see e.g. \cite{EtVa} for introduction of this condition).
These conditions mean in particular that $B$ and $C$ are diagonal on the corresponding spaces, respectively $V_1$ and $V_2$. 
In addition, $D$ has components on basis elements $E_{ij}\otimes E_{kl}$ of $\mf{gl}(n)\otimes \mf{gl}(n)$ only when the sets
$\{i,k\}$ and $\{j,l\}$ are equal (property ZW). In other words non-zero elements have identical \emph{unordered} multiplets 
of line and column indices.

For the consistency of the exchange relations the following set of coupled ``dynamical'' YB-equations is imposed.
\eb
A_{12} \ A_{13} \ A_{23} \ &=& \ A_{23} \ A_{13} \ A_{12} \label{dYBA}\\*
 D_{12}(\lambda+\gamma h_{3}) \ D_{13} \ D_{23}(\lambda+\gamma h_{1}) \ &=& \ D_{23} \ D_{13}(\lambda+\gamma h_2) 
\ D_{12} \label{dYBD}\\*
 D_{12} \ B_{13} \ B_{23}(\lambda+\gamma h_1) \ &=& \ B_{23} \ B_{13}(\lambda+\gamma h_2) \ D_{12} \label{dYBDB}\\*
 A_{12} \ C_{13} \ C_{23} \ &=& \ C_{23} \ C_{13} \ A_{12}(\lambda+\gamma h_3) \label{dYBAC}
\ee

The simplest example of this algebra is related to the elliptic scalar $\mf{gl}(n)$ Ruijsenaars-Schneider model and 
was first written in \cite{ACF}. 
We only write down its rational limit here.
\eb
A(\lambda)=1+\sum_{i\neq j}\frac{\gamma}{\lambda_{ij}} \left(E_{ii}-E_{ij}\right) \otimes \left(E_{jj}-E_{ji}\right) \label{RSexA}\\
B(\lambda)=C(\lambda)^{\pi}=1+\sum_{i\neq j}\frac{\gamma}{\lambda_{ij}-\gamma} E_{jj} \otimes \left(E_{ii}-E_{ij}\right)\\
D(\lambda)=1-\sum_{i\neq j}\frac{\gamma}{\lambda_{ij}} E_{ii} \otimes E_{jj} + \sum_{i\neq j}\frac{\gamma}{\lambda_{ij}} E_{ij} 
\otimes E_{ji}\label{RSexD}
\ee
where $E_{ij}$ is the elementary matrix whose entries are $\left(E_{ij}\right)_{kl}=\delta_{ik}\delta_{jl}$ and 
$\lambda_{ij}=\lambda_i-\lambda_j$.
These matrices verify the consistency conditions (\ref{dYBA})-(\ref{dYBAC}). 
A scalar representation of the exchange algebra defined with these structure matrices is then provided by:
\eb \ds
T(\lambda)=\sum_{ij} \frac{\prod_{a\neq i} (\lambda_{aj}+\tilde{\gamma})}{\prod_{a\neq j}\lambda_{aj}} E_{ij} \otimes \mathbf{1}
\label{RSexL}
\ee

The word ``scalar'' is used here in the sense that $T(\lambda)$ acts on a one-dimensional (trivial) quantum  space. 
The exchange relation (\ref{DQQA})
is just a c-number equality. Representation of (\ref{DQQA}) on non-trivial quantum spaces is provided in this context by 
the comodule structure in \cite{NAR}.  

Let us note here that the condition $AB=CD$ found in \cite{ACF} means in this context that the identity matrix 
is also a solution of (\ref{DQQA}). This is not a trivial statement; in fact it does not hold in general, and is not preserved
by fusion procedures.

\subsection{Fusion procedures and the ``dual'' algebra}

Let $A,B,C,D$ be solutions of the dynamical exchange relation. We will define their fusion by induction as follows. We omit the 
dependence on $\lambda$ and simplify the notations of the shifts as $(h_{(\ldots)})$; otherwise we use the notations 
introduced in section \ref{1fejez}, defining the multiple-index matrices by induction as:
\eb
&&A_{M\bar{N}'}=A_{1\bar{N}'}A_{M_0\bar{N}'}=A_{Mn'}A_{M{\bar{N}'}_{\phantom{,}0}} \non \\
&&B_{MN'}=B_{M1'}B_{MN'_0}=B_{1N'}\left[B_{M_0N'}(h_1)\right] \non \\
&&C_{MN'}=C_{1N'}C_{M_0N'}=C_{M1'}\left[C_{MN'_0}(h_{1'})\right] \non \\
&&D_{M\bar{N}'}=D_{1\bar{N}'}\left[D_{M_0\bar{N}'}(h_1)\right]=\left[D_{Mn'}(h_{(1',n'-1)})\right]D_{M{\bar{N}'}_{\phantom{,}0}} \non
\ee
where $\ds h_{(i,j)}:=\sum_{k=i}^j h_k$.
These fused structure matrices verify the fused dynamical YB-equations which are gathered together in the next proposition.
\begin{P1}
Let $A,B,C,D$ be solutions of the dynamical Yang-Baxter equations (\ref{dYBA})-(\ref{dYBAC}).Then the following fused dynamical Yang-Baxter equations hold:
\eb
A_{M\bar{N}'}A_{M\bar{L}''}A_{N'\bar{L}''} &=& A_{N'\bar{L}''}A_{M\bar{L}''}A_{M\bar{N}'} \non\\
A_{M\bar{N}'}C_{ML''}C_{N'L''} &=& C_{N'L''}C_{ML''}A_{M\bar{N}'}(h_{L''}) \non \\
D_{M\bar{N}'}(h_{L''}) D_{M\bar{L}''} D_{N'\bar{L}''}(h_{M}) &=& D_{N'\bar{L}''}D_{M\bar{L}''}(h_{N'})D_{M\bar{N}'} \non \\
D_{M\bar{N}'}B_{ML''}B_{N'L''}(h_{M}) &=& B_{N'L''}B_{ML''}(h_{N'})D_{M\bar{N}'} \non
\ee
\end{P1}
\bp by induction, using at crucial stages the zero weight properties. The fusion procedure respects the property ZW for $D$ and the
diagonality of $B$. It is also clear from the fusion procedure that the fused shift 
matrix $h_M$ is identified with $h_{(1,m)}$.
\end{proof}
\begin{T1}\label{theo7}
Let $T$ be a solution of the dynamical quadratic exchange relation 
\eb
A_{12}T_1B_{12}T_2(h_1)&=&T_2C_{12}T_1(h_2)D_{12}
\ee
then
\eb\label{dsolu}
 T_{M} = \prod^{\to}_{i \in M} \bigg( T_{i}(\sum_{\substack{k<i \\ k\in M}}h_k) \bigg( \prod^{\to}_{\substack{j>i\\ j\in M}} 
B_{ij}\bigg)\bigg)
\ee
verifies the fused dynamical exchange relation
\eb\label{dexch}
A_{M\bar{N}'}T_M B_{MN'}T_{N'}(h_M)=T_{N'}C_{MN'}T_M(h_{N'})D_{M\bar{N}'}
\ee
\end{T1}
\bp Similar to that of Theorem \ref{theo1} but the induction step uses the fact that $T_M=T_1 B_{1M_0} T_{M_0}(h_1)$ and
uses the fused dynamical YB-equations.  \end{proof}

The dual exchange relation and an associated fusion procedure are described in the next theorem.
\begin{T1}
Let $K$ be a solution of the dynamical quadratic exchange relation 
\eb\label{dexchDual}
(A_{12}^{-1})^{t_1t_2}K_1(B_{12}^{t_2})^{-1}K_2(h_1)&=&K_2(C_{12}^{t_1})^{-1}K_1(h_2)(D_{12}^{-1})^{t_1t_2}
\ee
then
\eb
 K_{M} = \prod^{\to}_{i \in M} \bigg( K_{i}(\sum_{\substack{k<i \\ k\in M}}h_k) \bigg( \prod^{\to}_{\substack{j>i\\ j\in M}} 
(B_{ij}^{t_j})^{-1}\bigg)\bigg)
\ee
verifies the fused dynamical exchange relation
\eb\label{dexchDual2}
\left(A_{M\bar{N}'}^{-1}\right)^{t_Mt_N}K_M \left(B_{MN'}^{t_{N'}}\right)^{-1}K_{N'}(h_M)&=&K_{N'}\left(C_{MN'}^{t_M}\right)^{-1}K_M(h_{N'})
\left(D_{M\bar{N}'}^{-1}\right)^{t_Mt_{N'}}
\ee
\end{T1}
\bp Similar to the nondynamical case.\end{proof} 

Note that the structure matrices of this dual relation are related to original ones in the same way as in the nondynamical case 
once we take into account the partial
zero weight property of $B$ and $C$ which implies diagonality on the corresponding spaces, respectively $V_1$ and $V_2$.

\subsection{Second fusion}

As in the nondynamical case, one can define another KDM-type fusion with the appropriate shifts. 
This fusion is characterized by the following exchange relation
\eb \label{dexch2}
A_{\bar{M}N'}T_M B_{M\bar{N}'} T_{N'}(h_M)&=&T_{N'} C_{\bar{M}N'} T_{M}(h_{N'}) D_{M\bar{N}'}
\ee

The analogy with the nondynamical case can be pushed further i.e.  there exists an object $L_M$ linking the fusions in Theorem \ref{theo7}
and \ref{theo9}. This allows us to use directly the proofs of Theorem \ref{theo3} and \ref{theo4}.

\begin{L1}
Let $T_M$ be a solution of the fused equation (\ref{dexch}). If $L_M$ verifies the following commutation rules
\eb
\label{kernel2}
L_MA_{M\bar{N}'}&=&A_{\bar{M}\bar{N}'} L_M\\
L_{N'}A_{M\bar{N}'}& =& A_{MN'} L_{N'} \non\\
L_{N'} B_{MN'}& =&B_{M\bar{N}'} L_{N'}(h_M) \non\\
L_{M} C_{MN'}&=&C_{\bar{M}N'} L_{M}(h_{N'})\non
\ee
then $L_MT_M$ is a solution of the exchange relation
\eb 
A_{\bar{M}N'}T_M B_{\bar{M}N'}T_{N'}(h_M)=T_{N'}C_{M\bar{N}'}T_M(h_{N'})D_{M\bar{N}'}
\ee
An example of such an $L_M$ is given by
\eb
\label{kernelEx2}
L_M=A_{12}\ldots A_{1m} A_{23}\ldots A_{2m} \ldots A_{m-1,m} = 
\prod^{\rightrightarrows}_{1\le i < j \le m} A_{ij}
\ee
\end{L1}
\bp Straightforward, using the dynamical YB-equations (\ref{dYBA})-(\ref{dYBDB}).\end{proof}

Now we state the dynamical versions of Theorem \ref{theo3} and \ref{theo4}.

\begin{T1}\label{theo9}
Let $T$ be a solution of the dynamical quadratic exchange relation
\eb
A_{12}T_1B_{12}T_2(h_1)&=&T_2C_{12}T_1(h_2)D_{12}
\ee
then
\eb\label{dsolu3}
T_M=\prod^{\to}_{i \in M}\bigg(\prod^{\to}_{\substack{j>i \\ j \in M}}A_{ij} T_i(\sum_{\substack{k<i\\k\in M}}h_k) 
\prod^{\gets}_{\substack{j>i \\ j \in M}} 
B_{ij}\bigg)
 \ee
verifies the fused dynamical exchange relation
\eb 
A_{\bar{M}N'}T_M B_{M\bar{N}'}T_{N'}(h_M)=T_{N'}C_{\bar{M}N'}T_M(h_{N'})D_{M\bar{N}'}
\ee
\end{T1}
\bp Reproduces the proof of Theorem \ref{theo3}, with suitable dynamical shifts.  \end{proof}

The dual exchange relation and an associated fusion procedure are described in the next theorem.
\begin{T1}
Let $K_M$ be a solution of the first fused exchange relation (\ref{dexchDual}) and $L_M$ be a solution of (\ref{kernel2}).
Then $\wt{K}_M=(L_M^{t_M})^{-1} K_M$ is a solution of the KDM-type dual fused exchange relation:
\eb \label{dFusEx3}
(A^{-1}_{\bar{M}N'})^{t_Mt_{N'}} \wt{K}_M  ((B^{t_M}_{M\bar{N}'})^{-1})^{t_{N'}}  \wt{K}_{N'}(h_{M})
 = \wt{K}_{N'} ((C_{\bar{M}N'}^{t_{N'}})^{-1})^{t_M}  \wt{K}_M(h_{N'})  (D_{M\bar{N}'}^{t_Mt_{N'}})^{-1}
\ee
\end{T1}
\bp Reproduces the proof of Theorem \ref{theo4}, with suitable dynamical shifts. \end{proof}

\subsection{Dressing.}

Solutions $T_M$ of the fused dynamical exchange relations also admit dressing procedures. However, because of the dynamical nature of the
exchange relations some of the equations that the dressings $Q_M$ and $S_M$ obey exhibit shifts, too. Specifically we have 
\begin{P1} \label{Prop6}
Let $T_M$ be a solution of the fused dynamical exchange relation. Then $Q_M T _M S_M$ is also a solution of the fused exchange relation
provided $Q_M$ and $S_M$ verify:

\eb
&\left[Q_M,A_{M\bar{N}'}\right]=\left[Q_{N'},A_{M\bar{N}'}\right]=0& \label{Dress3} \\
&Q_{N'}B_{MN'}=B_{MN'}Q_{N'}(h_{M}) \quad
Q_M C_{MN'}= C_{MN'}Q_M(h_{N'})& \non \\ 
&\left[S_{N'},C_{MN'}\right]=\left[S_M,B_{MN'}\right]=0& \\
&S_M(h_{N'}) D_{M\bar{N}'}=D_{M\bar{N}'}S_M \quad
S_{N'}D_{M\bar{N}'}=D_{M\bar{N}'}S_{N'}(h_M)& \non 
\ee
A particular solution of these constraints is given by:
\eb
&&Q_{M}=\check{A}_{12}\check{A}_{23}\ldots \check{A}_{m-1,m}\non\\
&&S_{M}=\check{D}_{12}\check{D}_{23}(h_{1})\ldots \check{D}_{m-1,m}(h_{(1,m-2)}),\non
\ee
\end{P1}
\bp By induction, similar to the non-dynamical dressings.\end{proof}
An interesting comparison can be drawn between this formula for $S_M$ and the formula used in \cite{ABB} to dress the quantum traces for
dynamical quantum groups.
The formula for $S_M$ is exactly the ``mirror image'' of the formula: $S^{ABB}_M=\check{R}_{12}(h_{(3,m)})\ldots\check{R}_{m,m-1}$.

\subsection{Three lemmas: dynamical and cyclic properties of D.}

Three easy technical lemmas are required to proceed with the construction.

\begin{L1}[Dynamical transposition]\label{DynTr}
Let $R(q)$ and $S(q)$ be two matrices with mutually commuting entries depending on a set of commuting coordinates $\ds \{q_k\}_{k=1}^n$. 
We then have:
\eb
\left(R(q)e^{\mc{D}}S(q)\right)^t =  \left[S^{SL}(q)\right]^t e^{\mc{D}} \left[R^{SC}(q)\right]^t
\ee
where $S^{SL}(q)_{ij}=e^{\partial_i} S(q)_{ij} e^{-\partial_i}=S_{ij}((q_1,\ldots,q_{\{i\}}+1,\ldots,q_n))$
(shift on line index) and $R^{SC}(q)_{ij}=e^{-\partial_j}R(q)_{ij}e^{\partial_j}$ (shift on column index).
\end{L1}
\bp We compare the $ij$-th entry on both sides using the fact that entries of $S^{SL}$ and 
$R^{SC}$ do not contain explicit shift quantities $e^{\partial}$ and therefore commute with each other.
If (in the case of $k$-tensor products) $''i''$ denotes a $k$-uple of indices $(i_1,\ldots i_k)$, the notation $q_{\{i\}}+1$
must be interpreted as $q_{i_1}+1, \ldots q_{i_k}+1$.
\end{proof}

{\bf Remark.} Later we will use this lemma in the special case when $R(q)$ is diagonal. This implies $R^{SC}(q)=
e^{-\mc{D}} R(q) e^{\mc{D}}$.

\begin{L1}[Matrix dynamical shift]\label{MatShift}
Let $D(q)$ be a matrix obeying the zero weight condition:
\eb
\left[D_{12},h \otimes \mathbf{1}+\mathbf{1} \otimes h\right]=0 \qquad \qquad (h \in \mathfrak{h})
\ee
Then the exponentials can be ``pushed through'' $D$, that is we have 
\eb
e^{-\D_1-\D_2} D_{12}=\bar{D}_{12}e^{-\D_1-\D_2}
\ee
where $\bar{D}_{12}=D_{12}^{-SL_{12}}$.
\end{L1}
\bp
What this lemma means is that one can write $e^{-\D_1-\D_2}D_{12}e^{\D_1+\D_2}$ in a \emph{matrix form} where
the exponentials of derivatives cancel out. The proof is straightforward
because the  zero weight condition implies the identification of incoming and outgoing indices of $D$. One then verifies easily the
equality of the two sides.
\end{proof}
\begin{L1}\label{Cycl}
Let $D(q)$ be a matrix obeying the zero weight condition:
\eb
\left[D_{12},h \otimes \mathbf{1}+\mathbf{1} \otimes h\right]=0 \qquad \qquad (h \in \mathfrak{h})
\ee
Then $D$ is cyclic with respect to the trace operation over $V_1\otimes V_2$ as follows:
\eb
Tr_{12} \left(D_{12} X_{12} D^{-1}_{12} e^{\mc{D}_1} e^{\mc{D}_2}\right)= Tr_{12}\left(X_{12} e^{\mc{D}_1} e^{\mc{D}_2}\right)
\ee
where $X$ is an arbitrary matrix the entries of which commute with the entries of $D$.
\end{L1}
\bp 
Consequence of the ZW property of $D$, which allows to reinterpret the matrix indices of $e^{\D_1+\D_2}$ 
as line instead of column indices of $D_{12}^{-1}$, allowing then to independently sum over the now decoupled column indices
of $D_{12}^{-1}$ with line indices of $D_{12}$ to altogether eliminate the matrix $D$ from the trace. Labels $1$ and $2$ 
formally denote here tensored auxiliary spaces.
\end{proof}

\subsection{Commuting hamiltonians.}

We can now state the fundamental result of this section.

\begin{T1}\label{theo11}
Let $\mathcal{T}_M$ be a solution of the fused dynamical exchange relations (\ref{dexch}). $\mc{T}_{M}$ acts on the tensor 
product of the auxiliary spaces labeled by $M$ and on the quantum space $V_q$.

Let $\mc{K}_M$ be a solution of the dual fused dynamical exchange relation (\ref{dexchDual2}). $\mc{K}_M$ acts on the tensor 
product of the auxiliary  spaces labeled by $M$ and on the quantum space $V_{q'}$.

The following operators
\eb
H_M=Tr_M \left(T_M e^{\mc{D}_M} (K^{SC}_M)^{t_M}\right)
\ee
constitute a family of mutually commuting quantum operators acting on $V_q \otimes V_{q'}$

\eb
\left[H_M,H_{N'}\right]=0
\ee

\end{T1}
\bp Similar to the preceding one, but extra care must be taken because of the shift operators that enter the expression. Using 
the dynamical transposition
lemma for $K_{N'}$ one has:
\eb
&&H_M H_{N'}= Tr \left( T_M e^{\mc{D}_M} (K^{SC}_M)^{ t_M}  T_{N'} e^{\mc{D}_{N'}} (K_{N'}^{SC})^{ t_{N'}}\right)= \\
&&Tr \left( T_M e^{\mc{D}_M} (K_M^{SC})^{ t_M}  T_{N'}^{t_{N'}} K_{N'} e^{\mc{D}_{N'}} \right) \non
\ee
since the invariance of the trace with respect to transposition is preserved in the dynamical case. 
One then writes:
\eb
&& Tr \left( T_M T_{N'}^{t_{N'}}(h_M)  \ e^{\mc{D}_M} (K_M^{SC})^{ t_M}  K_{N'} e^{\mc{D}_{N'}} \right)= \non \\
&& Tr\left( T_M T_{N'}^{t_{N'}}(h_M)
B_{MN'}^{t_{N'}} (B_{MN'}^{t_{N'}})^{-1} e^{\mc{D}_M} (K_M^{SC})^{ t_M} K_{N'} e^{\mc{D}_{N'}} \right)= \non \\
&& Tr\left( (T_M B_{MN'} T_{N'}(h_M))^{t_M t_{N'}} ((B_{MN'}^{t_{N'}})^{-1} e^{\mc{D}_M} (K_M^{SC})^{ t_M})^{t_M} K_{N'} e^{\mc{D}_{N'}} 
\right) \non
\ee
In the last equality the identification $ T_{N'}^{t_{N'}}(h_M) B_{MN'}^{t_{N'}}= (B_{MN'} T_{N'}(h_M))^{t_{N'}}$ uses the 
zero-weight condition
$\left[B_{MN'},h_{M}\right]=0$. Using once again the dynamical transposition lemma and the zero-weight condition on $B$ which 
guarantees $e^{\mc{D}_M}((B_{MN'}^{t_{N'}})^{-1})^{SCt_M}=(B_{MN'}^{t_{N'}})^{-1} e^{\mc{D}_M}$ as commented above, one gets:
\eb
Tr\left( (T_M B_{MN'} T_{N'}(h_M))^{t_M t_{N'}}A_{MN'}^{t_Mt_{N'}} 
(A_{MN'}^{t_Mt_{N'}})^{-1} K_M (B_{MN'}^{t_{N'}})^{-1} 
e^{\mc{D}_M} K_{N'} e^{\mc{D}_{N'}} \right)= \non \\
 Tr\left( (A_{MN'}T_M B_{MN'} T_{N'}(h_M))^{t_M t_{N'}}(A_{MN'}^{t_Mt_{N'}})^{-1} K_M (B_{MN'}^{t_{N'}})^{-1} 
  K_{N'}(h_M) e^{\mc{D}_M}e^{\mc{D}_{N'}} \right)= \non
\ee
One here identifies the direct and dual exchange relation. to yield:.
\eb
Tr \left( ( T_{N'} C_{MN'} T_M(h_{N'}) D_{MN'})^{t_Mt_{N'}}  K_{N'} (C_{MN'}^{t_M})^{-1} K_M(h_{N'}) 
(D_{MN'}^{t_Mt_{N'}})^{-1} e^{\mc{D}_M}e^{\mc{D}_{N'}} \right)= \non \\
Tr \left(  D_{MN'}^{t_Mt_{N'}}( T_{N'} C_{MN'} T_M(h_{N'}))^{t_Mt_{N'}}  K_{N'} (C_{MN'}^{t_M})^{-1} K_M(h_{N'}) 
(D_{MN'}^{t_Mt_{N'}})^{-1} e^{\mc{D}_M}e^{\mc{D}_{N'}} \right) \non 
\ee
Here Lemma \ref{Cycl} is at work.
\eb
&&Tr \left( ( T_{N'} C_{MN'} T_M(h_{N'}))^{t_Mt_{N'}}  K_{N'} (C_{MN'}^{t_M})^{-1} K_M(h_{N'}) 
e^{\mc{D}_M}e^{\mc{D}_{N'}} \right)= \non \\
&&Tr \left( T_{N'} (C_{MN'} T_M(h_{N'}))^{t_M}  (K_{N'} (C_{MN'}^{t_M})^{-1} e^{\mc{D}_{N'}})^{t_{N'}} K_M 
e^{\mc{D}_M} \right)= \non \\ 
&&Tr \left( T_{N'} T_M^{t_M}(h_{N'}) C_{MN'}^{t_M}  (C_{MN'}^{t_M})^{-1} e^{\mc{D}_{N'}} K_{N'}^{SC t_{N'}} 
K_M 
e^{\mc{D}_M} \right) \non
\ee
Once again we have used the dynamical transposition lemma and the partial weight zero property of $C_{MN'}$  

\eb
&&Tr \left( T_{N'} T_M^{t_M}(h_{N'})  e^{\mc{D}_{N'}}  (K_{N'}^{SC})^{ t_{N'}}    K_M 
e^{\mc{D}_M} \right)= \non \\
&& Tr \left( T_{N'} e^{\mc{D}_{N'}} (K_{N'}^{SC})^{ t_{N'}} T_{M}^{t_{M}}  K_{M} e^{\mc{D}_M}  \right)= 
Tr\left(  T_{N'} e^{\mc{D}_{N'}} K_{N'} T_{M} e^{\mc{D}_M} (K_M^{SC})^{ t_M} \right)\non
\ee
\end{proof}
Without the dressing described by Proposition \ref{Prop6} the traces constructed in (\ref{dsolu}) decouple just as in the nondynamical case.
Indeed we have
\begin{P1}
Operators built from the solution (\ref{dsolu}) decouple as $$Tr_M(T_Me^{\mc{D}_M} (K_M^{SC}(^{ t_M})= Tr(T e^{\mc{D}} 
(K^{SC})^{ t})^{\# M}$$
\end{P1}
\bp We will prove the proposition for $M$ with two elements. The statement remains valid for higher powers by induction.
We also need to put the trace under a more amenable form. In fact, $Tr(T_M e^{\mc{D}_M}(K_M^{SC})^{ t_M})=Tr(T_M^{t_M}K_M e^{\mc{D}_M})$.
by virtue of Lemma \ref{DynTr}.
\eb
&&Tr(\left[T_1 B_{12} T_2(h_1)\right]^{t_1t_2} K_1 (B_{12}^{t_2})^{-1} K_2(h_1)e^{\mc{D}_1+\mc{D}_2})= \non \\
&&Tr(T_1 \left[ B_{12}  T_2(h_1)\right]^{t_2} \left[ K_1 e^{\mc{D}_1} ((B^{t_2}_{12})^{-1})^{SC_1}\right]^{t_1}K_2 e^{\mc{D}_2})= \non
\ee
where ${(\phantom{0})}^{SC_1}$ means $(\phantom{0})^{SC}$ operation applied on the first space.
\eb
&&Tr(T_1 T_2^{t_2}(h_1)B_{12}^{t_2}  (B_{12}^{t_2})^{-1} e^{\mc{D}_1} K_1^{SC t_1}K_2e^{\mc{D}_2})= \non
Tr(T_1  e^{\mc{D}_1} T_2^{t_2} K_1^{SC t_1} K_2e^{\mc{D}_2})=\\
&&Tr(T_1 e^{\mc{D}_1} K_1^{SC t_1} T_2^{t_2} K_2 e^{\mc{D}_2}) = Tr(T e^{\mc{D}}K^{SC t})^2 \non 
\ee
\end{proof}

Of course, the three comments made after Proposition \ref{decoup} in the non-dynamical case remain valid, 
although we do not know yet of explicit examples for Mezincescu-Nepomechie procedure in a dynamical context.

\section{The fully dynamical algebra} 

The third type of quadratic algebra considered here is the extension to general structure matrices $A,B,C,D$ of the ``boundary
dynamical algebra'' (BDA) considered in \cite{MPLA, Gould}. Fusion and trace formulas were defined in \cite{MPLA} 
for the particular case of BDA
where $A=D=R(u_1-u_2)$, $B=C=R(u_1+u_2)$, $R$ being the IRF $\mathbb{Z}_n$ $R$-matrix. The most general ``fully dynamical'' 
(denomination to be justified presently) exchange algebra reads:
\eb\label{ExRel}
A_{12}(\lambda) T_{1}(\l+\g h_2) B_{12}(\l) T_{2}(\l+\g h_1) = T_{2}(\l+\g h_1)  C_{12}(\l)  T_{1}(\l+\g h_2)  D_{12}(\l)
\ee

Once again we assume $dim \ \mf{h} = dim \ V$ \cite{EtVa}. 
The following conditions are imposed on the structure matrices ($R=A$, $B$, $C$ or $D$)
\bigskip
Unitarity
\eb
R_{12}(u_1,u_2;\l) R_{21}(u_2,u_1;\l)=\mathbf{1}
\ee

Zero weight property
\eb \label{zw}
[ h \otimes \mathbf{1}+\mathbf{1} \otimes h, R_{12}(u_1,u_2;\l)]=0 \qquad (h \in \mf{h}) 
\ee

$R$ then verifies
the same ZW property as in the semi-dynamical case.
By contrast with the previous case all four matrices in (\ref{ExRel}) exhibit a dynamical shift and all four structure
matrices have $(1+2)$- zero weight, hence the denomination ``fully dynamical''.
In some specific examples \cite{Gould,Fe} the structure matrices also obey the dynamical zero weight property:
\eb \label{dzw}
[\mc{D} \otimes \mathbf{1}+\mathbf{1} \otimes \mc{D}, R_{12}(u_1,u_2;\l)]=0 
\ee

Structure matrices all obey Gervais-Neveu-Felder type equations.

\eb \label{YBfd}\non
A_{12}(\l) A_{13}(\l+\g h_2) A_{23}(\l) & = & A_{23}(\l+\g h_1) A_{13}(\l) A_{12}(\l+\g h_3) \\ \non
A_{12}(\l) C_{13}(\l+\g h_2) C_{23}(\l) & = & C_{23}(\l+\g h_1) C_{13}(\l) A_{12}(\l+\g h_3) \\ \non
D_{12}(\l+\g h_3) D_{13}(\l) D_{23}(\l+\g h_1) & = & D_{23}(\l) D_{13}(\l+\g h_2) D_{12}(\l) \\ 
D_{12}(\l+\g h_3) B_{13}(\l) B_{23}(\l+\g h_1) & = & B_{23}(\l) B_{13}(\l+\g h_2) D_{12}(\l)  
\ee

If the dynamical zero weight property is verified then all equations can be rewritten under the more familiar 'alternating shift' form 
$$R_{12}(\l-\g h_3)R_{13}(\l+\g h_2)R_{23}(\l-\g h_1)=R_{23}(\l+\g h_1)R_{13}(\l-\g h_2)R_{12}(\l+\g h_3)$$ 

As in the previous situation, these equations ensure the compatibility of the algebra in the following sense.
Let us take the left hand side of exchange relation (\ref{ExRel}), embed it in a triple tensor product and shift it on the third space. 
Then let
us multiply it  with $B_{13}(u_1,u_3;\l) B_{23}(u_2,u_3;\l +\g h_1) T_3(u_3;\l+\g h_1+\g h_2)$. One can reverse the order of the $T$'s
in two different ways which yield the same result if equations (\ref{YBfd}) are obeyed.


\subsection{Fusion procedure and the ``dual'' algebra}

The fusion of the structure matrices is again defined by induction as follows:

\eb
&&A_{M\bar{N}'}=A_{1\bar{N}'}(h_{(2,m)})A_{M_0\bar{N}'}=A_{Mn'}A_{M\bar{N}'_0}(h_{n'}) \non \\
&&B_{MN'}=B_{1N'}B_{M_0N'}(h_1)=B_{M1'}(h_{(2',n')})B_{MN'_0} \non \\
&&C_{MN'}=C_{1N'}(h_{(2,m)})C_{M_0N'}=C_{M1'}C_{MN'_0}(h_1) \non \\
&&D_{M\bar{N}'}=D_{1\bar{N}'}D_{M_0\bar{N}'}(h_1)=D_{Mn'}(h_{(1',n'-1)})D_{M\bar{N}'_0} \non
\ee

These fused matrices verify the corresponding fused YB-equations and the ZW property.

\begin{P1}
Let $A,B,C,D$ be solutions of the dynamical Yang-Baxter equations (\ref{YBfd}).Then the following fused dynamical Yang-Baxter 
equations hold:
\eb
A_{M\bar{N}'}A_{M\bar{L}''}(h_{N'})A_{N'\bar{L}''} &=& A_{N'\bar{L}''}(h_{M})A_{M\bar{L}''}A_{M\bar{N}'}(h_{L''}) \non\\
A_{M\bar{N}'}C_{ML''}(h_{N'})C_{N'L''} &=& C_{N'L''}(h_{M})C_{ML''}A_{M\bar{N}'}(h_{L''}) \non \\
D_{M\bar{N}'}(h_{L''}) D_{M\bar{L}''} D_{N'\bar{L}''}(h_{M}) &=& D_{N'\bar{L}''}D_{M\bar{L}''}(h_{N'})D_{M\bar{N}'} \non \\
D_{M\bar{N}'}(h_{L''})B_{ML''}B_{N'L''}(h_{M}) &=& B_{N'L''}B_{ML''}(h_{N'})D_{M\bar{N}'} \non
\ee
\end{P1} 
\bp
 Straightforward by induction.
\end{proof}

Note that the dynamical zero weight property 
does not survive fusion, but algebraic zero weight does.
In this sense this dynamical zero weight property is not relevant for the construction of commuting traces, 
and is not (generically) a feature
of the universal algebra. We will from now on disregard it. In addition, we will concentrate here on the most relevant features of quantum
trace building, ignoring for instance the possibility of a ``second fusion''.

\begin{T1}
Let $T$ be a solution of the dynamical quadratic exchange relation 
\eb
A_{12}T_1(h_2)B_{12}T_2(h_1)&=&T_2C_{12}T_1(h_2)D_{12}
\ee
then
\eb
 T_{M} = \prod^{\to}_{i \in M} \bigg( T_{i}(\sum_{\substack{k\neq i \\ k\in M}}h_k) \bigg( \prod^{\to}_{\substack{j>i\\ j\in M}} 
B_{ij}(\sum_{\substack{k<i\\k \in M}}h_k+\sum_{\substack{k>j \\ k\in M}}h_k) \bigg)\bigg)
\ee
verifies the fused dynamical exchange relation
\eb
A_{M\bar{N}'}T_M(h_{N'})B_{MN'}T_{N'}(h_M)=T_{N'}(h_M)C_{MN'}T_M(h_{N'})D_{M\bar{N}'}
\ee
\end{T1}
\bp Similar to that of Theorem \ref{theo1} but the induction step uses the fact that $T_M=T_1(h_{M_0}) B_{1M_0} T_{M_0}(h_1)$ and
uses the fused dynamical YB-equations.  \end{proof}

The dual exchange relation and the associated fusion procedure are described in the next theorem.
\begin{T1}
Let $K$ be a solution of the dynamical quadratic exchange relation
\eb\non
A^d_{12}(\l) K_1(\l+\g h_2) B^d_{12}(\l)) 
K_2(\l+\g h_1) =  \non
K_2(\l+\g h_1) C^d_{12}(\l) K_1(\l+\g h_2) 
D^d_{12}(\l)
\ee
where 
\eb
&A^{d}_{12}=((A_{12}^{-SL_{12}})^{-1})^{-SC_{12}t_{12}}\quad B^{d}_{12}=(((B_{12}^{-SL_{12}})^{-SC_2t_2})^{-1})^{SL_1t_1}& \non \\
&C^d_{12}=(((C_{12}^{-SL_{12}})^{-SC_1t_1})^{-1})^{SL_2t_2} \quad D^d_{12}=((D_{12}^{-SL_{12}})^{-1})^{SL_{12}t_{12}} \non &
\ee
then
\eb\label{dsolu2}
 K_{M} = \prod^{\to}_{i \in M} \bigg( K_{i}(\sum_{\substack{k\neq i \\ k\in M}}h_k) \bigg( \prod^{\to}_{\substack{j>i\\ j\in M}} 
B^d_{ij}(\sum_{\substack{k<i\\k \in M}}h_k+\sum_{\substack{k>j \\ k\in M}}h_k) \bigg)\bigg)
\ee
verifies the fused dynamical exchange relation
\eb
A^d_{M\bar{N}'}K_M(h_{N'})B^d_{MN'}K_{N'}(h_M)=K_{N'}(h_M)C^d_{MN'}K_M(h_{N'})D^d_{M\bar{N}'}
\ee
\end{T1}
\bp Straightforward once one has established that the fused dual structure matrix is equal to the dual of the fused structure matrix and
that the YB-equations obeyed by the dual structure matrices derive from the equations (\ref{YBfd}).
\end{proof}
\subsection{Dressing}

\begin{P1} 
Let $T_M$ be a solution of the fused fully dynamical exchange relation. 
Then $Q_M T _M S_M$ is also a solution of the fused exchange relation
provided $Q_M$ and $S_M$ verify:

\eb
&Q_M A_{M\bar{N}'}=A_{M\bar{N}'}Q_M(h_{N'}) \quad Q_{N'}(h_M)A_{M\bar{N}'}=A_{M\bar{N}'} Q_{N'}&  \\
&Q_{N'}B_{MN'}=B_{MN'}Q_{N'}(h_{M}) \quad
Q_M C_{MN'}= C_{MN'}Q_M(h_{N'})& \non \\ 
&S_{N'}(h_M)C_{MN'}=C_{MN'}S_{N'} \quad S_M(h_{N'}) B_{MN'}=B_{MN'} S_M& \\
&S_M(h_{N'}) D_{M\bar{N}'}=D_{M\bar{N}'}S_M \quad
S_{N'}D_{M\bar{N}'}=D_{M\bar{N}'}S_{N'}(h_M)& \non 
\ee
A particular solution of these constraints is given by:
\eb
&&Q_{M}=\check{A}_{12}(h_{(3,m)})\check{A}_{23}(h_{(4,m)})\ldots \check{A}_{m-1,m}\non\\
&&S_{M}=\check{D}_{12}\check{D}_{23}(h_{1})\ldots \check{D}_{m-1,m}(h_{(1,m-2)}),\non
\ee
\end{P1}
\bp By induction.
\end{proof}

\subsection{Commuting traces}

We use the following properties inferred from lemma \ref{MatShift}.
\eb\non
e^{-\D_1-\D_2} A_{12}=A_{12}^{-SL_{12}} e^{-\D_1-\D_2}=\bar{A}_{12}e^{-\D_1-\D_2} \\ \non
e^{-\D_2} A_{12} e^{\D_1} = e^{\D_1} \bar{A}_{12} e^{-\D_2} 
\ee
and their transposed variants:
\eb
e^{\D_1} (\bar{A}_{12}^{-SC_2t_2})^{-1} e^{\D_2}= e^{\D_2} (A_{12}^{-SL_2t_2})^{-1} e^{\D_1}
\ee
and so on.
Since these relations are immediately derived from the ZW property on the structure matrices, they remain valid for fused 
structure matrices, 
too, since the fusion respects the zero weight property as opposed to the dynamical zero weight property (cf. remark above). 
In this case labels $1$ and $2$ formally denote tensored auxiliary spaces.

\begin{T1}\label{theo14}
Let $\mathcal{T}_M$ be a solution of the fused dynamical exchange relations (\ref{dexch}). $\mc{T}_{M}$ acts on the tensor 
product of the auxiliary spaces labeled by $M$ and on the quantum space $V_q$.

Let $\mc{K}_M$ be a solution of the dual fused dynamical exchange relation (\ref{dexch2}). $\mc{K}_M$ acts on the tensor 
product of the auxiliary  spaces labeled by $M$ and on the quantum space $V_{q'}$.

The following operators
\eb
H_M=Tr_M e^{-\D_M} T_M e^{\D_M} K_M^{SCt_M}
\ee
constitute a family of commuting operators acting on $V_q \otimes V_{q'}$

\eb
\left[H_M,H_{N'}\right]=0
\ee

\end{T1}
\bp
It is worth to give a detailed description of the proof as in theorem \ref{theo11} since the occurence 
of derivative objects $\sim e^{\D_M}$ 
considerably modifies it in comparison to the standard Sklyanin-type proof for non-dynamical
algebras. Once again the dynamical transposition lemma plays a essential role.
\eb \non
&&H_MH_{N'}=Tr \  e^{-\D_M} T_M e^{\D_M} K_M^{SC_Mt_M}  \ \ e^{-\D_{N'}} T_{N'} e^{\D_{N'}} K_{N'}^{SC_{N'}t_{N'}} =\\ \non
&&Tr \ e^{-\D_M} T_M e^{\D_M} K_M^{SC_Mt_M}  \  \left( e^{-\D_{N'}} T_{N'}\right)^{t_{N'}}\left( e^{\D_{N'}} K_{N'}^{SC_{N'}t_{N'}}\right)^{t_{N'}} =\\ \non
&&Tr \  e^{-\D_M} T_M e^{\D_M} K_M^{SC_Mt_M} \ T_{N'}^{-SL_{N'}t_{N'}} e^{-\D_{N'}} K_{N'} e^{\D_{N'}} =\\ \non
&&Tr \ e^{-\D_M} T_M e^{\D_M} T_{N'}^{-SL_{N'}t_{N'}} \  \ K_M^{SC_Mt_M} e^{-\D_{N'}} K_{N'} e^{\D_{N'}} =\\ \non
&&Tr \ e^{-\D_M} \left[e^{-\D_{N'}} A_{MN'}^{-1} T_{N'}(h_M) C_{MN'} T_M(h_{N'}) D_{MN'} e^{\D_M}\right]^{t_{N'}} \times \non \\
&&\qquad (\bar{B}_{MN'}^{-SC_{N'}t_{N'}})^{-1} e^{\D_{N'}}  \non
K_M^{SC_Mt_M} e^{-\D_{N'}} K_{N'} e^{\D_{N'}} = \\ \non
&&Tr \left[e^{-\D_M-\D_{N'}} A_{MN'}^{-1} T_{N'}(h_M) C_{MN'} T_M(h_{N'}) D_{MN'} \right]^{t_{MN'}} \times \non \\
&&\qquad \left[e^{\D_M} (\bar{B}_{MN'}^{-SC_{N'}t_{N'}})^{-1} e^{\D_{N'}} K_M^{SC_Mt_M}\right]^{t_M} e^{-\D_{N'}} K_{N'} e^{\D_{N'}} \non
\ee
Pushing exponentials through $B$.
\eb \non
&&Tr \left[e^{-\D_M-\D_{N'}} A_{MN'}^{-1} T_{N'}(h_M) C_{MN'} T_M(h_{N'}) D_{MN'} \right]^{t_{MN'}} e^{\D_{N'}} \times \non \\
&& \qquad\left[(B_{MN'}^{-SL_{N'}t_{N'}})^{-1} e^{\D_M} K_M^{SC_Mt_M}\right]^{t_M} e^{-\D_{N'}} K_{N'} e^{\D_{N'}} = \non \\ \non
&&Tr \left[\bar{A}_{MN'}^{-1}e^{-\D_M-\D_{N'}} T_{N'}(h_M) C_{MN'} T_M(h_{N'}) D_{MN'} \right]^{t_{MN'}} e^{\D_{N'}} \times \non \\
&& \qquad K_M e^{\D_M} 
((B_{MN'}^{-SL_{N'}t_{N'}})^{-1})^{SC_Mt_M} e^{-\D_{N'}} K_{N'} e^{\D_{N'}} \non
\ee
Using zero weight of $A$ and $B$ transposed.
\eb \non
&&Tr \left[T_{N'}(h_M) C_{MN'} T_M(h_{N'}) D_{MN'} \right]^{-SL_{MN'}t_{MN'} }e^{-\D_M-\D_{N'}} \times\\ 
\non 
&& \qquad (\bar{A}_{MN'}^{-1})^{-SC_{MN'}t_{MN'}} e^{\D_{N'}} K_M 
e^{-\D_{N'}} ((\bar{B}_{MN'}^{-SC_{N'}t_{N'}})^{-1})^{SL_1t_M} e^{\D_M} K_{N'} e^{-\D_M} e^{\D_M+\D_{N'}} =\\ \non
&& Tr \left[T_{N'}(h_M) C_{MN'} T_M(h_{N'}) D_{MN'} \right]^{-SL_{MN'}t_{MN'} }e^{-\D_M-\D_{N'}} \times \\ \non
&&\qquad \{(\bar{A}_{MN'}^{-1})^{-SC_{MN'}t_{MN'}}K_M(h_{N'})
((B_{MN'}^{-SC_{N'}t_{N'}})^{-1})^{SL_1t_M} K_{N'}(h_M)\} e^{\D_M+\D_{N'}} =\\ \non
&& Tr \left[T_{N'}(h_M) C_{MN'} T_M(h_{N'}) D_{MN'} \right]^{-SL_{MN'}t_{MN'} }e^{-\D_M-\D_{N'}} \times \\
&& \qquad K_{N'}(h_M) ((\bar{C}_{MN'}^{-SC_Mt_M})^{-1})^{SL_{N'}t_{N'}} \non
K_M(h_{N'}) (\bar{D}_{MN'}^{-1})^{SL_{MN'}t_{MN'}} e^{\D_M+\D_{N'}} =\\ \non
&& Tr \left[e^{-\D_M-\D_{N'}}T_{N'}(h_M) C_{MN'} T_M(h_{N'})D_{MN'} \right] \times \non \\
&& \qquad \left[K_{N'}(h_M) ((\bar{C}_{MN'}^{-SC_Mt_M})^{-1})^{SL_{N'}t_{N'}}
K_M(h_{N'}) e^{\D_M+\D_{N'}}(D_{MN'}^{-1})^{SC_{MN'}t_{MN'}}\right]^{t_{MN'}} = \non\\ \non
&& Tr e^{-\D_M-\D_{N'}}T_{N'}(h_M) C_{MN'} T_M(h_{N'})D_{MN'} D_{MN'}^{-1}  e^{\D_M+\D_{N'}} \times \\ \non
&& \qquad \left[K_{N'}(h_M) ((\bar{C}_{MN'}^{-SC_Mt_M})^{-1})^{SL_{N'}t_{N'}} K_M(h_{N'})\right]^{SL_{MN'}t_{MN'}} =\\ \non
&& Tr e^{-\D_{N'}} T_{N'} \left[e^{-\D_M} C_{MN'} e^{\D_{N'}}T_M\right]^{t_M} e^{\D_M-\D_{N'}} \times \\ \non
&&\qquad \left[K_{N'} e^{-\D_M} ((\bar{C}_{MN'}^{-SC_Mt_M})^{-1})^{SL_{N'}t_{N'}} 
e^{\D_{N'}} \right]^{t_{N'}} K_M e^{\D_M} =\\ \non 
&&Tr e^{-\D_{N'}} T_{N'} e^{\D_{N'}} \left[\bar{C}_{MN'} e^{-\D_M} T_M\right]^{t_M} e^{\D_M-\D_{N'}} \times \\ \non 
&& \qquad \left[K_{N'} e^{\D_{N'}} (C_{MN'}^{-SL_1t_M})^{-1})^{SC_{N'}t_{N'}}
\right]^{t_{N'}} e^{-\D_M} K_M e^{\D_M} =\\ \non
&& Tr e^{-\D_{N'}} T_{N'} e^{\D_{N'}} T_M^{-SL_1t_M} e^{-\D_M} \bar{C}_{MN'}^{-SC_Mt_M} e^{\D_M-\D_{N'}}(C_{MN'}^{-SL_1t_M})^{-1} e^{\D_{N'}} K_{N'}^{SC_{N'}t_{N'}}
e^{-\D_M} K_M e^{\D_M}
\ee
Using zero weight of $C$.
\eb \non
&&Tr e^{-\D_{N'}} T_{N'} e^{\D_{N'}} T_M^{-SL_1t_M} e^{-\D_{N'}} C_{MN'}^{-SL_1t_M}(C_{MN'}^{-SL_1t_M})^{-1} e^{\D_{N'}} K_{N'}^{SC_{N'}t_{N'}}
e^{-\D_M} K_M e^{\D_M} =\\ \non
&&Tr e^{-\D_{N'}} T_{N'} e^{\D_{N'}} T_M^{-SL_1t_M} K_{N'}^{SC_{N'}t_{N'}} e^{-\D_M} K_M e^{\D_M} =\\ \non
&&Tr e^{-\D_{N'}} T_{N'} e^{\D_{N'}} K_{N'}^{SC_{N'}t_{N'}} T_M^{-SL_1t_M} e^{-\D_M} K_M e^{\D_M} =\\ \non
&&Tr e^{-\D_{N'}} T_{N'} e^{\D_{N'}} K_{N'}^{SC_{N'}t_{N'}} \left[T_M^{-SL_1t_M} e^{-\D_M}\right]^{t_M} \left[K_M e^{\D_M}\right]^{t_M} =\\ \non
&&Tr e^{-\D_{N'}} T_{N'} e^{\D_{N'}y} K_{N'}^{SC_{N'}t_{N'}} e^{-\D_M} T_M e^{\D_M} K_M^{SC_Mt_M}=H_{N'}H_M
\ee
\end{proof}
\section{Conclusion}

We have now defined fusion and trace procedures in view of obtaining commuting hamiltonians of ``quantum trace type'', 
for the non-dynamical general quadratic algebra (\ref{QQA}),  for the semi-dynamical quadratic algebra (\ref{DQQA}) and
for the fully dynamical quadratic algebra (\ref{ExRel}). Our
immediate interest is now to apply this procedure to some particularly interesting examples of such
quadratic algebras, the most relevant being at this time the scalar Ruijsenaars-Schneider quantum Lax formulation (semi-dynamical type)
\cite{NA}.

Note in this respect that previous application of an order-one trace formulation (i.e without auxiliary space 
tensor products) to the specific case of ``boundary dynamical
$\mathfrak{sl}(2)$ algebras'' considered in \cite{Ra} yielded models described in \cite{Gould} as generalizations of the Gaudin models.
Positions of the sites were associated with values of the spectral parameters (in a spin-chain type construction), not with the
dynamical variable itself whose interpretation is unclear.

As already emphasized, our elucidation of tensor product structure for quadratic algebras is also very important in formulating
generalizations of the Mezincescu-Nepomechie fusion procedure in general open spin chains \cite{Annecy}.

Our constructions moreover also shed light on some characteristic properties of the quadratic algebra. The building of 
commuting traces requires first of all the introduction
of a dual exchange relation. It seems possible that this notion reflects the existence of anti-automorphisms of the
underlying hypothetical algebra structure, of which the transposition and crossing-relations used in the
non-dynamical cases (see \cite{AD}) would be realizations.

The explicit formulation of consistent fusion relations should also help in understanding the meaning of
quantum algebra (QA) structures and characterizing in particular their coalgebra properties. As pointed out, the DKM-type fusions do stem
in at least one case from a universal structure \cite{DKM}, and so does the fusion for boundary dynamical algebra (case 
when $A,B,C,D$ stem from one single dynamical $R$-matrix \cite{KM}). 
Regarding the semi-dynamical QA 
it was already known \cite{NAR} that one could extend the quantum space on which entries of $T$ act, by auxiliary spaces of
$A$ and $B$ or $C$ and $D$ matrices, thereby obtaining spin-chain like construction of a monodromy matrix (comodule structure).
We have now defined the complementary procedure, extending the auxiliary space by a ``fusion'' procedure. This yields
the full ``coproduct'' or rather comodule structure of the DQA (\ref{DQQA}).

\vspace{5mm}
 {\bf Acknowledgments}:
AD is supported by the TMR Network ``EUCLID''; ``Integrable models and applications: from strings to condensed matter'', contract
number HPRN-CT-2002-00325.

 


\begin{thebibliography}{99}

\bibitem{cherednik}  I. ~Cherednik, Theor. Math. Phys. \textbf{61} (1984), 977; E. ~K. ~Sklyanin, J. Phys.   
\textbf{A21} (1988), 2375; D. ~Fioravanti, M. ~Rossi, J. Phys \textbf{A34} (2001), 567; M. ~Mintchev, E. ~Ragoucy,
P. ~Sorba, \emph{hep-th/0303187}; A.~Kundu, Mod. Phys. Lett. \textbf{A10} (1995), 2955; 
L. ~Hlavaty, Journ. Math. Phys. \textbf{35} (1994)
, 2560; S. ~Majid, Journ. Math. Phys \textbf{32} (1991), 3246.

\bibitem{AD} J. Avan, A. Doikou, Commuting quantum traces: the case of reflection algebras, J. Phys \textbf{A 36} (2003), p. 1; \emph{math.QA/0305424}

\bibitem{FrMa91a} L.~Freidel, J.~M.~Maillet, Quadratic algebras and integrable systems, Phys. Lett \textbf{B 262} (1991), p. 268.

\bibitem{Fad} L. Faddeev, Integrable models in (1+1)-dimensional quantum field theory, in: Les Houches 1982 ed.:J. B. Zuber and R. Stora
pp. 561-608


\bibitem{ABB} J. Avan, O. Babelon, E. Billey, The Gervais-Neveu-Felder equation and the quantum Calogero-Moser systems, Comm. Math. Phys. 
\textbf{178} (1996), p. 281; \emph{hep-th/9505091}

\bibitem{Maillet} J.M. Maillet, Lax equations and quantum groups, Phys. Lett. \textbf{B 245} (1990), p. 480.

\bibitem{Ma85} J.~M.~Maillet, Phys. Lett \textbf{B 162} (1985), p. 137.

\bibitem{MeNe}L.~Mezincescu, R.~I.~Nepomechie, Fusion procedure for open spin chains, J. Phys. \textbf{A 25} (1992) p. 2533.

\bibitem{Annecy}  
D. Arnaudon, J. Avan , N.~Cramp\'e, 
  A.~Doikou, L. Frappat, {E}. Ragoucy, 
General boundary conditions for the $sl({\cal N})$ 
 and super $sl({\cal M}/{\cal N})$ open spin chains, J. Stat. Mech.: Theor. Exp. (JSTAT) P08 (2004) P08005, \emph{math-ph/0406021}

\bibitem{KiRe} A. N. Kirillov, N. Yu. Reshetikhin, Exact solution of the integrable XXZ Heisenberg model, J.Phys. \textbf{A 20} (1987) 
p. 1565; P. P. Kulish, N. Yu. Reshetikhin and E.K. Sklyanin, Lett. Math. Phys. \textbf{5} (1981) p. 393.

\bibitem{NAR} Z. Nagy, J. Avan and G. Rollet, Construction of dynamical quadratic algebras, Lett. 
Math. Phys. \textbf{67} (2004) p. 1; \emph{math.QA/0307026}

\bibitem{Ra} 
Heng Fan, Bo-Yu Hou, Kang-Jie Shi, Representation of the boundary elliptic quantum group $\mc{B}E_{\tau, \eta}(sl_2)$ and the
Bethe ansatz, Nucl. Phys. {\bf B 496} (1997) p. 551-570; 

\bibitem{MPLA} Heng Fan, Bo-You Hou, Guang-Liang Li, Kang-Jie Shi, Integrable $A^{(1)}_{n-1}$ IRF model with reflecting boundary
conditions Mod. Phys. Lett. A \textbf{26} (1997) pp. 1929-1942.

\bibitem{KM} P. P. Kulish, A. I. Mudrov: Dynamical reflection equation, \emph{math.QA/0405556}

\bibitem{DKM} J.~Donin, A.~I.~Mudrov, Reflection equation, twist and equivariant quantization, Isr. J. Math. \textbf{136} (2003), 
p. 11., \emph{math.QA/0204295};  J.~Donin, P.~P.~Kulish and A.~I.~Mudrov, 
On universal solution to reflection equation,  Lett.Math.Phys., 
\textbf{63} (2003) p.179; \emph{math.QA/0210242}

\bibitem{Sk} P.~P.~Kulish, E.~K.~Sklyanin: Algebraic structure related to the reflection equation, J. Phys. \textbf{A 25} 1992)
p. 5963.

\bibitem{FrMa91b}  L.~Freidel, J.~M.~Maillet On classical and quantum integrable field theorioes associated to Kac-Moody current algebras, Phys. Lett \textbf{B 263} (1991), p. 403.


\bibitem{doikou} A. Doikou, J. Phys. {\bf A33} (2000) 8797.

\bibitem{mac} J.B. McGuire, J. Math. Phys. {\bf 5}
(1964) 622.

\bibitem{yang} C.N. Yang, Rev. Lett. {\bf 19} (1967) 1312.

\bibitem{baxter}
  R.J. Baxter,  Ann. Phys. \textbf{70} (1972) 193; J. Stat. Phys.
  \textbf{8} (1973) 25; {\it Exactly solved models in statistical
    mechanics} (Academic Press, 1982)

\bibitem{korepin}
  V.E. Korepin, Theor. Math. Phys. \textbf{76} (1980) 165; V.E.
  Korepin, G. Izergin and N.M. Bogoliubov, {\it Quantum inverse
    scattering method, correlation functions and algebraic Bethe
    Ansatz} (Cambridge University Press, 1993).


\bibitem{zamo} A.B. Zamolodchikov, Al.B. Zamolodchikov, Ann.Phys. {\bf120} (1979) 253.

\bibitem{faddeev}  L.D. Faddeev and L.A. Takhtajan, J.Sov.Math. {\bf 24} (1984)
241;\\ L.D. Faddeev and L.A. Takhtajan, Phys.Lett. {\bf 85A} (1981)
375.

\bibitem{EtVa}P. Etingof, A. Varchenko, Solution of the quantum dynamical Yang-Baxter equation and dynamical quantum groups, 
Comm. Math. Phys \textbf{196} (1998) p. 591, \emph{q-alg/9708015}

\bibitem{ACF} G.E.Arutyunov, L.O. Chekhov, S.A. Frolov, R-matrix quantization of the elliptic Ruijsenaars-Schneider model, Comm. Math. Phys.
{\bf 192} (1998), pp. 405-432, {\it q-alg/9612032}; G. E. Arutyunov, S. A. Frolov, Quantum dynamical R-matrices and quantum Frobenius group
Comm. Math. Phys. {\bf 191} (1998), pp. 15-29, {\it q-alg/9610009}

\bibitem{Gould}
Mark D. Gould, Yao-Zhong Zhang, Shao-You Zao, Elliptic Gaudin models and elliptic KZ equations, Nucl.Phys. {\bf B 630} (2002) p.492-508, 
\emph{nlin.SI/0110038}; 

\bibitem{Fe} G.Felder , Proc. ICM Z\"urich \emph{hep-th/9407154} (1994), 1247; Proc. ICMP Paris (1994), 211. 


\bibitem{NA} Z. Nagy, J. Avan: Spin chains from dynamical quadratic algebras, in preparation. 

\end{thebibliography}
\end{document}